\documentclass[12pt]{amsart}
\usepackage[backref=page,colorlinks,bookmarksopen=true]{hyperref}
\usepackage{amssymb}
\usepackage{latexsym}
\usepackage{amsmath}

\newcommand{\eenmatrix}[4]{\setlength{\unitlength}{1ex}
\raisebox{-1.5ex}[0ex][3.5ex]{\parbox{7.5ex}{\begin{picture}(5,5)(0,0)
\put(1.8,1.9){\line(0,1){4}}
\put(1.8,5.9){\line(1,0){4}}
\put(5.8,5.9){\line(0,-1){4}}
\put(5.8,1.9){\line(-1,0){4}}
\put(2.3,6.2){\parbox[b]{3ex}{\centering $#1$}}
\put(6.2,3.4){\parbox[b]{25ex}{$#2$ }}
\put(-3,3.4){\parbox[b]{1.5ex}{  \raggedleft $#3$  }}
\put(-1,0){\parbox[t]{10ex}{\centering $#4$}}
\end{picture}

}}}
\newcommand{\eenmatrixbase}[4]{\setlength{\unitlength}{1ex}
\raisebox{-1.5ex}[0ex][3.5ex]{\parbox{7.5ex}{\begin{picture}(5,5)(0,0)
\put(1.8,1.9){\line(0,1){4}}
\put(1.8,5.9){\line(1,0){4}}
\put(5.8,5.9){\line(0,-1){4}}
\put(5.8,1.9){\line(-1,0){4}}
\put(2.3,6.2){\parbox[b]{3ex}{\centering $#1$}}
\put(6.2,3.4){\parbox[b]{8ex}{$#2$ }}
\put(0,3.4){\parbox[b]{1.5ex}{  \raggedleft $#3$  }}
\put(2.3,0){\parbox[t]{3ex}{\centering $#4$}}
\end{picture}
}}}
\newcommand{\eenmatrixbis}[4]{\setlength{\unitlength}{1ex}
\raisebox{-1.5ex}[0ex][3.5ex]{\parbox{7.5ex}{\begin{picture}(5,5)(0,0)
\put(1.8,1.9){\line(0,1){4}}
\put(1.8,5.9){\line(1,0){4}}
\put(5.8,5.9){\line(0,-1){4}}
\put(5.8,1.9){\line(-1,0){4}}
\put(2.3,6.2){\parbox[b]{3ex}{\centering $#1$}}
\put(6.2,3.4){\parbox[b]{12ex}{$#2$ }}
\put(-3,3.4){\parbox[b]{1.5ex}{  \raggedleft $#3$  }}
\put(-1,0){\parbox[t]{10ex}{\centering $#4$}}
\end{picture}

}}}
\newcommand{\eenmatrixfevrier}[4]{\setlength{\unitlength}{1ex}
\raisebox{-1.5ex}[0ex][3.5ex]{\parbox{7.5ex}{\begin{picture}(5,5)(0,0)
\put(1.8,1.9){\line(0,1){4}}
\put(1.8,5.9){\line(1,0){4}}
\put(5.8,5.9){\line(0,-1){4}}
\put(5.8,1.9){\line(-1,0){4}}
\put(2.3,6.2){\parbox[b]{3ex}{\centering $#1$}}
\put(6.2,3.4){\parbox[b]{25ex}{$#2$ }}
\put(0,3.4){\parbox[b]{1.5ex}{  \raggedleft $#3$  }}
\put(2.3,0){\parbox[t]{3ex}{\centering $#4$}}
\end{picture}
}}}

\begin{document}

\title{ Relative  matched pairs of finite groups from   depth two inclusions of von Neumann algebras to quantum groupoids }
\author{Jean-Michel Vallin}
\address{UMR CNRS 6628 Universit\'eŽ d'Orl\'eŽans \\  Institut de
Math\'eŽmatiques de Chevaleret \\Plateau 7D, 175
rue du Chevaleret 75013 Paris
\\e-mail: jmva @ math.jussieu.fr}

\date{Preliminary version of the 03/27/07}
\subjclass{(2000) 46L37, 20L05, 20G42, 20 99}
\keywords{ Subfactors, quantum groupoids, relative matched pairs.}

\markboth{Jean-Michel Vallin}{ Relative  matched pairs of finite groups and quantum groupoids with a non abelian  basis }

\begin{abstract} In this work we     give  a generalization    of  matched pairs of (finite) groups  to describe  a general class of      depth two inclusions of  factor von Neumann algebras and the     C*-quantum groupoids associated with, using  double groupoids.
\end{abstract}

\maketitle

\newpage
\section{Introduction}
\def\I{ I_{\mathcal H,\mathcal K}}
\def\k{\gamma^K}
\def\h{\chi^H}
\newenvironment{dm}{\hspace*{0,15in} {\it Proof:}}{\qed}

The first examples of C*-quantum groupoids were discovered by the theoretical physicists B\"ohm, Szlachanyi and Nill  (\cite{BoSz}, \cite{BoSzNi}, they called them weak Hopf $C^*$-algebras. In the infinite dimensional  framework and for a very large class of 
depth two inclusions of von Neumann algebras, we proved  with M.Enock \cite{EV}  the existence of such objects, giving  a non commutative geometric   interpretation of their basic construction.
D.Nikshych  and  L.Vainerman, using general inclusions of depth two subfactors of 
type $II_1$ with finite index  $M_0 \subset M_1$, gave, thanks to a specific bracket,    weak Hopf $C^*$-algebra structures to  the relative commutants $M'_0 \cap M_2$ and $M'_1\cap M_3$,  \cite{NV1}, M.C. David in \cite{D} gave also a refinement of this property.

In this article, using  a generalization of  matched pairs of groups inside the theory of groups, we prove the existence of a large class of inclusions of von Neumann algebras  in the conditions of Vainerman and Nykshych.  As we have no reference,  we shall prove, thanks to Dietmar Bisch works (\cite{B1}, \cite{BH}),  the following property which may be well known by specialists:  
\vskip 0.4cm
{\bf {Proposition}}
{\it Let $H$ and $K$ be two finite subgroups of a group   $G$ such that $G = HK = \{hk /h \in H , k\in K\}$, let $R$ be the type $II_1$ factor, if $G$ acts properly and outerly on $R$ and if one denotes by $R^H$ (resp. $R \rtimes K$) the fixed point algebra (resp. crossed product) of the induced action of $H$(resp. $K$) then  the inclusion  $R^H \subset R \rtimes K$ is a  depth two inclusion of subfactors of 
type $II_1$ with finite index.}

Thanks to  the   canonical bracket, we  describe simply the C*-quantum groupoids associated with this inclusion,     using a  double groupoid   structure in the sense of Ehresmann \cite{EH} and a pair of C*-quantum groupoid structure associated with  $H,K$. These structures are very close to  examples due to  N. Andruskiewitsch and S.Natale ( \cite{AN1}  and \cite{AN2}),  so we  give  a precise relation, they ask for,  in the  introduction of \cite{AN2},  between  matched pairs and their  work. Hence,  our aim is to prove the following theorem:

{\bf {Theorem }}
{\it Let $H$ and $K$ be two finite subgroups of a group   $G$ such that $G = HK = \{hk /h \in H , k\in K\}$, then the   pair of C*-quantum groupoids in duality associated  with the inclusion $R^H \subset R \rtimes K$ is isomorphic to   the double groupoid's pair.  }

To be clear, the article do not follow the above order (the one of our study). 

In the second paragraph we recall   the definition  of a C*-quantum groupoid and we explain the commutative and symmetric examples.

In the third chapter, we define relative  matched pairs of groups, we associate  to them two canonical double crossed product  and two    double groupoid structures in duality; this allows us to define explicit C*-quantum groupoid in duality for this double crossed products. One can see, that  the technics we  use with  sets of representatives, generalize, in a certain way, ideas developped independently in  \cite{B},  for the particular case of inclusions of groups. Then we prove these pairs  give also depth two inclusions  of von Neumann algebras   and  so  C*-quantum groupoids in duality in an other way. 
The fourth chapter proves that these  two structures are isomorphic.

So a natural extension of this article will be the generalization of these constructions in the direction of Lesieur's locally compact groupoids \cite{L}. An other will be a characterization of these objects in terms of cleft extensions in the spirit of S.Vaes and L.Vainerman \cite{VV}.

I want to thank a lot L.Vainerman having suggested me this  study,   D.Bisch for some explanations about the type of  inclusion we deal with, and also S.Baaj, M.C.David and M.Enock for the numerous discussions we had.

\section{ C*-quantum groupoids}
\label{decadix}

Let us recall the definition of a C*-quantum groupoid (or a weak Hopf 
$C^\star$ -algebra), one can also see a more synthetic  approach ( \cite{Val1}\cite{Val2}....)  using a generalization of Baaj and Skandalis's multiplicative unitaries (\cite{BS}Êand \cite{BBS}). We shall also recall the definition  of an action of such a C*-quantum groupoid due to  D.Nikshych and L.Vainerman (\cite{NV1}) (alternative definitions can be found in {\cite[Chap.3]{Val2} or {\cite{E2}, which are suitable for a future generalisation to Lesieur's measured quantum groupoid theory \cite{L}).

\subsection{C*-quantum groupoids }
\subsubsection{\bf{Definition}}(G.B\"ohm, K.Szlach\'anyi, F.Nill) 
(\cite{BoSz}, \cite{BoSzNi})\label{Bohm}

{\it A weak Hopf  $C^*$-algebra is a collection $(A, \Gamma, \kappa,
\epsilon)$ where: $A$ is a finite-dimensional $C^*$-algebra (or von Neumann algebra), $\Gamma: 
A \mapsto A \otimes A$ is a generalized coproduct, which means that: 
$(\Gamma\otimes i)\Gamma = (i \otimes \Gamma)\Gamma$, $\kappa$ is an
antipode on $A$, i.e., a linear map  from $A$ to $A$ such 
that $(\kappa \circ *)^2 = i$ (where $*$ is the involution on $A$), 
$\kappa(xy) = \kappa(y)\kappa(x)$ for every  $x,y$ in $A$ with
$(\kappa \otimes \kappa) \Gamma = \varsigma\Gamma \kappa$ (where 
$\varsigma$ is the usual flip on $A\otimes A$). 

We suppose also that 
$(m(\kappa \otimes i)\otimes i)(\Gamma \otimes i)\Gamma(x) = 
(1 \otimes x)\Gamma(1)$ (where $m$ is the multiplication of
tensors, i.e., $m(a \otimes b) = ab$), and that $\epsilon$ is a   
counit, i.e., a positive linear form on $A$ such that
$(\epsilon \otimes i)\Gamma = (i\otimes \epsilon)\Gamma = i$, and 
for every  $x,y$ in $A$: $(\epsilon \otimes\epsilon) ((x \otimes 1) 
\Gamma(1)(1 \otimes y)) = \epsilon(xy)$. }

\subsubsection{\bf{Results}} (cf. \cite{NV1},\cite{BoSzNi})
\label{denuit}
{\it If $(A, \Gamma, \kappa, \epsilon)$ is a weak Hopf  $C^*$-algebra,
then the sets 
$ A_t =\{ x\in A/\Gamma(x) = \Gamma(1)(x\otimes 1) = 
(x \otimes 1)\Gamma(1), \} $ and 
$A_s= \{ x \in A / \Gamma(x) = 
\Gamma(1)(1\otimes x) = (1\otimes x)\Gamma(1)\} $
are commuting sub $C^*$-algebras of $A$ and $ \kappa(A_t) = A_s$; one calls them  respectively target and 
source  Cartan subalgebra of $(A, \Gamma, \kappa, \epsilon)$ or   simply { \bf  basis} of $(A, \Gamma, \kappa, \epsilon)$.}

In fact, we shall here deal only with the special case of C*-quantum groupoids for which $\kappa$ is involutive, namely weak Kac algebras.

\subsection{ The commutative and symmetric examples}  \label{metz}

\vskip 1cm
Let's recall that a groupoid $\mathcal G$ is a small category the  morphisms of which are all invertible. In  all what follows, $\mathcal G$ is finite. Let $\mathcal G^0$ be the   set of objects,  one can identify  $\mathcal G^0$ to a subset of the morphisms. So a (finite)  groupoid can also be viewed as a set $\mathcal G$ together with a, not everywhere defined,  multiplication for which there is a set of unities $\mathcal G^0$, two maps, source denoted by $s$ and target by $t$, from $ \mathcal G$  to $\mathcal G^0$ so that the product $xy$ of two elements $x,y \in \mathcal G $ exists  if and only if $s(x) = t(y)$;  every element $x \in \mathcal G$ has a unique inverse $x^{-1}$, and one has $x(yz) = (xy)z$ whenever both members make sense. We refer to \cite{R} for  the fondamental structures and notations for groupoids.

Let's denote $H = l^2(\mathcal G)$, with the usual notations. One can define two C*-quantum groupoids in duality acting on $H = l^2(\mathcal G)$, namely $(C(\mathcal G), \Gamma_{\mathcal G}, \kappa_{\mathcal G},\epsilon_{\mathcal G})$, the commutative example, and $(\mathcal R(\mathcal G), \hat \Gamma_{\mathcal G}, \hat \kappa_{\mathcal G}, \hat \epsilon_{\mathcal G}) $, the symmetric example,  where: $ C(\mathcal G)$ is  the
commutative involutive algebras of complex valued functions on
$\mathcal G$, $  \mathcal R(\mathcal G) $ = $ \{  \sum_{x \in {\mathcal G}} a_x \rho (x )\} $ is the right regular algebra of $\mathcal G$  and $\rho(x) $ is the  partial isometry given by the formula $(\rho(x)\xi)(t) = \xi(tx)$ if $x \in \mathcal G^{s(t)}$ and $= 0$ otherwise. 
The two  $ *$-quantum groupoids structures on $S$ and $\hat S$ are given by:
 \vskip 0,5cm
\begin{itemize}

\item \   \    \   \    \    \  {\bf Coproducts}:

$\Gamma_{\mathcal G}(f)(x,y) = f(xy)$ if  x,y  are composables and $f \in C(\mathcal G)$.

 \   \    \   \    \    \   \   \    \   \    \    \  \   \    \   \     \   \    \    \    $=0$   \   \     \   \    \    \   \  \  otherwise

$\hat \Gamma_{\mathcal G}(\rho(s)) =  \rho(s) \otimes \rho(s) $

\item \   \    \   \    \    \  {\bf Antipodes:}

$\kappa_{\mathcal G}(f)(x) = f(x^{-1}),$  \  \  $Ê\hat \kappa_{\mathcal G}(\rho(s)) =  \rho(s^{-1}) =  \rho(s)^* $

\item \   \    \   \    \    \  {\bf Counities:}

$\epsilon_{\mathcal G}(f) = \sum_{u\in \mathcal G^0}f(u)$,   \  \  \  $ \hat \epsilon_{\mathcal G}(\rho(s)) = 1$

\end{itemize}

To define the symmetric example, one also could consider the left regular representation of $\mathcal G$: using the $*$-algebra $ \mathcal L(\mathcal G)= \{  \sum_{s \in {\mathcal G}} a_s \lambda(s) \} $  ( the left regular algebra of $\mathcal G$), where $\lambda(s) $ is the  partial isometry given by the formula $(\lambda(s)\xi)(t) = \xi(s^{-1}t)$ if $t \in \mathcal G^{r(s)}$ and $= 0$ otherwise.

\subsection{Action of a C*-quantum groupoid }
\subsubsection{\bf{Definition}} 
(\cite{NV1} chap 2)
\label{beamer}
{\it A (left) action of a C*-quantum groupoid $A$ on a von Neuman algebra $M$ is any linear weakly continous map:
$ A \otimes M \to M : a \otimes m \mapsto a\triangleright m$ defining on  $M$  a left $A$-module structure and such that :

i) $a \triangleright(xy) = (a_{(1)}\triangleright x)(a_{(2)}\triangleright y)$, where $\Gamma(a) = a_{(1)}\otimes a_{(2)}$ with Sweedler notations,

ii) $(a \triangleright x)^* =  \kappa(a) \triangleright (x^*)$ 

iii) $a \triangleright 1 =  \epsilon^t(a) \triangleright 1$, and $a \triangleright 1 = 0$ iff  $\epsilon^t(a) = 0$ (where $ \epsilon^t(a) =  \epsilon(1_{(1)}a)1_{(2)}$)  }
\vskip 0.1cm
Thanks to the antipode, $M$  also becomes a right module ($m.a= \kappa(a) \triangleright m$). When such an action is given one can define the fixed point algebra as the von Neumann  $M^A = \{ m \in M | \forall a \in A, a \triangleright m = \epsilon^t(a) \triangleright m \}$. Also the von Neumann crossed product $ M \rtimes A$ (see \cite{NV1} chap 2) is the $\mathbb C$-vector space $M \underset {A_t} \otimes A$ where $M$ is the right $A_t$ module $M$ obtained via the right multiplication by the elements $a_t \triangleright 1$ ($a_t \in A_t$) and $A$ is the left $A_t$ module obtained by the left $A_t$ multiplication. So elements of $M \underset {A_t} \otimes A$ are the classes $ [m \otimes a]$ with the identification: $m(a_t \triangleright 1) \otimes a \equiv m \otimes a_ta$. The multiplication and the involution of $M \underset {A_t} \otimes A$ are given for any $m,m' \in M$ and $a,a' \in A$ by the formulas:
$$ [m \otimes a][m' \otimes a'] = [m(a_{(1)} \triangleright m') \otimes a_{(2)}a']$$
$$[m \otimes a]^* =  [(a^*_{(1)} \triangleright m^*) \otimes a_{(2)} ]$$

\subsection{  C*-quantum groupoids in action}
\label{ses}
We suppose known the Jones s' tower theory (see \cite{GHJ}) and in this section we recall the tight relation between  C*-quantum groupoids  actions and  depth two inclusions of  type $II_1$-factors with finite index $\lambda^{-1} = [M:N]$, (see \cite{NV1} for full detail). 

Let $M_0 \subset M_1$ be such an inclusion.

 Let $M_0 \subset M_1 \subset M_2 = <M_1,e_1> \subset  M_3 = <M_2,e_2> \subset...$ be the Jones construction. 
 
 Let $A = M'_0 \cap M_2$, $B= M'_1 \cap M_3$, and let $\tau$ be  the Markov trace of the inclusion; if one defines for every $a \in A $ and $b \in B$ the bracket:
  $$<a,b> = d \lambda^{-2}\tau(ae_2e_1b)$$
 where $d = dim(M'_0 \cap M_1)$, this defines a non degenerate duality between $A$ and $B$, the   following theorem is true:
 
 \subsubsection{\bf{Theorem}}
 \label{tenir}
 (\cite{NV1}{th 4.17 and chap 6), \cite{D} 3.8.4 }

 {\it  For every   $b \in B$, let's define $\Gamma_B(b)$,$\epsilon_B(b)$,$\kappa_B(b)$, such that for any $a,a' \in A$ one has:
$$ <aa',b> = <a \otimes a',\Gamma_B(b)>$$
$$\epsilon_B(b) = <1,b>  $$
$$ <a,\kappa_B(b) = \overline{<a^*,b^*>}$$

and let's suppose that $\Gamma_B$ is multiplicative, then:

i)  $(B,\Gamma_B ,\kappa_B,\epsilon_B)$   is a  C*-quantum groupoid, 

ii) the map $\triangleright: B \otimes M_2 \mapsto M_2:  b \triangleright x = \lambda^{-1}E_{M_2}(bxe_2)$ 
defines a left action of $B$ on $M_2$, $M_1$ is  the fixed point subalgebra $M_2^B$ and the map $\theta : [x \otimes b] \mapsto xb$ is an isomorphism between  the crossed product  $M_2 \rtimes B$ and $M_3$.

iii) using  the above bracket leads to define a dual C*-quantum groupoid on $A$ and an action of $B$ on $A$ which is the   restriction of the action $\triangleright$ to $A$. }

\subsubsection{\bf{Corollary}}
\label{jncg}
{\it The map $[a \otimes b] \mapsto ab $ is an isomorphism between the crossed product $A\rtimes B$ and $M'_0 \cap M_3$.}
\newline
\begin{dm} Clearly the map $[a \otimes b] \mapsto ab $ is an isomorphism, its image is included in $M'_0 \cap M_3$; conversely, due to the depth two condition, the Jones  projection $e_2$ is also the first Jones projection for  the derived tower: $M'_0 \cap M_1 \subset M'_0 \cap M_2 \subset M'_0 \cap M_3$, hence $M'_0 \cap M_3$ is generated by $M'_0 \cap M_2$ and $e_2$ so it is included in (and equal to)   the image of the isomorphism.
\end{dm}

\vskip 0.5cm

\section{Relative  matched  pairs of groups }
 Let's now explain what we mean by a relative  matched pair.
  
  \vskip 0,2cm
  
  \subsection{Relative  matched pairs of groups}
\label{simplifier}
  \subsubsection{\bf{Definition}}
\label{Aisey}
{\it  Let $G$ be a  group, two any subgroups $H,K$ of $G$ are said to be a relative  matched pair  if and only if $G = HK = \{hk/ h\in H , k\in K\}$.  }
    \vskip 0.6cm
\subsubsection{\bf{Remark and notations}}
\label{lorient}
A relative  matched pair $H,K$ is a  matched pair if and only if $H \cap K = \{e\}$ where $e$ is the unit   of $G$. For the sake of simplicity, let's denote $ {\bf S = H \cap K}$. Of course one can construct a lot of examples of relative  matched pairs:  if $H$ is any subgroup of $G$ then $H,G$ is a relative  matched pair different from a  matched pair if $H \not = \{e\}$. Let's give a nice machinery to obtain examples:
    \vskip 0.6cm
\subsubsection{\bf{Lemma (Frattini Argument, th 1.11.8 of \cite{G})}}
\label{superbe}
{\it  Let $G$ be a finite group, let $N$ be a normal subgroup of $G$, let $P$ be a Sylow p-subgroup of $N$, then $G = NN_G(P)$, where $N_G(P)$ is the normalizer of $P$ in $G$.}

    \vskip 0.6cm
\subsubsection{\bf{Lemma and notations}}
\label{muet}
{\it    Let $H,K$ be a relative  matched pair in $G$, for any $(h,k)$ in $H \times K$ let's denote  $p_1(hk)  =  hS$ and $p_2(hk) = Sk$, this defines   two   maps            $p_1: G \to \ \  H_{\mid S}$ and $p_2: G \to \ \ _{S\mid}K$ such that:

i)  for any $g$ in $G$ and any $(h,k)$ in $p_1(g)\times p_2(g)$ there exists a unique $(h',k')$ in $p_1(g)\times p_2(g)$ verifying $g = hk' = h'k$

ii) for any $g,g'$ in $G$, one has: $p_1(g) = p_1(g')$ (resp.$p_2(g) = p_2(g')$) if and only if there exists $ k \in K$ (resp.$h \in H$) such that $g' = gk$ (resp. $g' = hg$). }
\vskip 0,3cm
\begin{dm}
For any $g \in G$, there exists $ (h,k)$ in $H \times K$ such that $g =hk$ and if $(h',k')$ in $H \times K$ verifies  $hk = h'k'$, this exactly means there exists $t \in S$ such that $h' =ht$ and $k' = t^{-1}k$, the existence of $p_1$ and $p_2$  and the first assertion of the lemma follow. For any $g,g'$ in $G$, let $(h',k') \in H\times K$ such that $g= h'k'$, so  if $p_1(g') = p_1(g)$ then  there exists  $k \in K$ such that $g' = h'k$, so $g' = (h'h^{-1})g$, conversely, if there exists $k \in K$ such that $g' = gk$ then $g' = h'k'k$, hence $h' \in p_1(g')$ so $p_1(g) = p_1(g')$.
\end{dm}
\vskip 0.5cm

\subsubsection{\bf{Remark and notation}} 
\label{demoralise}

 As $K,H$ is also a relative  matched pair, there exists also two   maps   $p'_1: G \to K_{\mid S}$ and $p'_2: G \to \  Ö  _{ S\mid}H$,  defined for any $(h,k)$ in $H \times K$ by $p'_1(kh) = kS$ and $p'_2(kh) = Sh$

  \vskip 1.5cm
\subsection{\bf{Double crossed products associated with relative  matched pairs}}
\subsubsection{\bf{Lemma}}
\label{lesieur}
{\it    For any $(h,k)$ in $H\times K$,  let $ k \triangleright   hS$ (resp.$ h \triangleright' k S$)  be equal to $  p_1(kh)$  (resp.$   p'_1(hk)$), this defines   a   left action of  group $K$  on the set  $H_{\mid S}$ and a   left action of  group $H$  on the set  $K_{\mid S}$.}
 \vskip 0.6cm

\begin{dm}
Easy.
\end{dm}
\vskip 0.6cm
 \subsubsection{\bf{Notations}}
\label{noter}
Till the end of this article we shall consider an exhaustive family  $I = \{ k_1, k_2,...\}$  (resp. $J = \{h_1,h_2,...\} $) of representatives of $K_{\mid S}$ (resp. $H_{\mid S}$). We shall now extend the actions above:

 \subsubsection{\bf{Proposition}}
\label{marteau}
{ \it Let $ \underset I \triangleright'$  be defined  for any $h \in H$, $k_i \in I$ and $s \in S$ by:
$$ h  \underset I \triangleright' k_is = k_js  \  \  \  \   where  \  \  \{k_j \} = I \cap  h   \triangleright'  k_iS $$
then $ \underset I \triangleright'$ is an action of $H$ on $K$ and for all $k \in K$, $s \in S$, one has: $ h  \underset I \triangleright' (ks)  = (h  \underset I \triangleright' k)s$ and $ h  \underset I \triangleright' s  = s$, if $I'$ is an other exhaustive family, the actions $ \underset I \triangleright'$ and $ \underset {I'}\triangleright'$ are conjugate.}
\vskip 0.6cm
\begin{dm}
Of course, $ \underset I \triangleright'$ is well defined, for any $k$ in $K$, $e \underset I \triangleright'  k = k$, and for any $h'$ in $H$, one has: $h'\underset I \triangleright'( h  \underset I \triangleright' k_is) = h'\underset I \triangleright'(k_js) = k_ls$ where : 
$$ 
 \{k_l \} 
 = I \cap  h'   \triangleright'  k_jS = I \cap h'   \triangleright' (h   \triangleright'  k_iS)  
 =  I \cap (h'h)  \triangleright'  k_i S
$$
 So $h'\underset I \triangleright'( h  \underset I \triangleright' k_is = (h'h) \underset I \triangleright' k_is$ and $\underset I \triangleright'$ is an action. For all $k \in K$, $s \in S$, let $  k_i   \in I  $ and $\sigma \in S$, such that $k = k_i\sigma$, one  has: $ h  \underset I \triangleright' (ks) = h  \underset I \triangleright' (k_i\sigma s) = k_j\sigma s$ where $  \{k_j \} = I \cap  h   \triangleright'  k_iS$, hence $ h  \underset I \triangleright' (ks) = (k_j \sigma)s= (h  \underset I \triangleright' k)s$. If $s_0 = I \cap  eS$  then $ h  \underset I \triangleright' s = h  \underset I \triangleright' s_0({s_0 }^{-1}s) = s_0({s_0 }^{-1}s) = s$.

Now let $I^1,I^2$ be two exhaustive families, for $i=1,2$, and any $c \in K_{\mid S}$, let $k^i_c = I^i \cap C$. We can define a permutation $\phi$ of $K$ such that for all  $c \in K_{\mid S}$ and $s \in S$, one has : $\phi(k^1_cs) = k^2_cs$,  then for all $h \in H$: $ h  \underset {I^2 } \triangleright' \phi(k^1_cs) = h \underset {I^2} \triangleright' k^2_cs = k^2_{h \triangleright'c} s = \phi(k^1_{h \triangleright'c} s) =  \phi(h \underset {I^1}\triangleright' k^1_cs)$ and $\phi$ realizes a conjugation between the two actions.
\end{dm}
\vskip 1cm
\subsubsection{{\bf Notations}} 
\label{andros}i) By reversing $H$ and $K$, one can  also extend the action $ \triangleright$ to an action $\underset J \triangleright$, of $K$ on $H$, with the same properties as in proposition \ref{marteau}.

ii) For any $h \in H$ and $k \in K$ let $ h \underset I \triangleleft' k $ and  $ k \underset J \triangleleft h $ be the unique element in $K$ and $H$ respectively such that: $hk = (h \underset I \triangleright' k )(h \underset I \triangleleft' k )$ and $kh  = (k \underset J \triangleright  h )(k  \underset J \triangleright h) $.

\vskip 0.3cm

One must keep in mind that in general $\underset I \triangleleft' $ and $\underset J \triangleleft$ are not ( right) actions.

 Let's define  a double crossed product with a relative  matched pair. In fact one can extend  the  action $\triangleright'  $ to the crossed product $ C(K) \underset {\rho} \rtimes  S$ of $K$ by  the right action of $S$, which is the $*$-algebra generated by a group of unitaries  $(\rho(s))_{s \in S}$ and a partition of the unity $(\chi_k)_{k \in K}$ with the commutation relations: $\rho(s)\chi_k = \chi_{ks^{-1}}\rho(s)$.  
\vskip 1cm

\subsubsection{\bf{Proposition}}
\label{roy}
{\it   For any $h \in H$, $k\in K$ and $s \in S$ , let's define:

$$ \sigma^I_h(\rho(s)\chi_k) = \rho(s)\chi_{h \underset I  \triangleright' k}   $$

then $(\sigma^I_h)_{h \in H}$, is an action of $H$ on the crossed product $  C(K) \underset {\rho} \rtimes S $, if $I'$ is an other exhaustive family, the actions $(\sigma^I_h)_{h \in H}$ and $(\sigma^{I'}_h)_{h \in H}$ are conjugate. }

\begin{dm} 

Obviously for any $h \in H$,  $\sigma^I_h$ is a well defined linear endomorphism on  $  C(K) \underset {\rho} \rtimes S $. For any $k\in K$ and $s \in S$, due to 1), one has : 
 $ \sigma^I_h(\rho(s)\chi_k\rho(s)^*) = \sigma^I_h(\chi_{ks^{-1}}) =  \chi_{h \underset I  \triangleright' (ks^{-1})} = \chi_{(h \underset I  \triangleright' k)s^{-1}} =   \rho(s)\sigma^I_h(\chi_k)\rho(s)^*$, the proposition follows.

\end{dm}
\vskip 0.5cm
So the double crossed product $ (C(K) \underset {\rho} \rtimes  S)\underset  {\sigma^I} \rtimes H$ can be viewed as the  $*$-algebra generated by the families $(\rho(s))_{s \in S}$, $(\chi_k)_{k \in K}$ and a group of unitaries $(V_h)_{h \in H}$ with the additional commutation relations: $\rho(s)V_h= V_h \rho(s)$ and $V_h\chi_k = \chi_{h \underset I  \triangleright' k}V_h$. Also one can extend   the action $\underset J \triangleright$ to an action $(\sigma^J_k)_{k \in K}$ of $K$ on the crossed product $  C(H) \underset {\rho} \rtimes S $, with the same properties.

\vskip 0,6cm
\subsection{{\bf Double groupoid structures and quantum groupoids structures associated with relative  matched pairs}}

\subsubsection{{\bf Definition}} Let $\mathcal T$  be the set : $ \{ \eenmatrixbase{h}{k}{k' }{h' } / \hskip 0.2cm   h,h' \in H , k,k' \in K   \hskip 0.1cm   \hskip 0.1cm  hk = k'h' \}$
 and let $\mathcal T'$ be   the set : $ \{ \eenmatrixbase{k}{h}{h' }{k' } / \hskip 0.2cm   h,h' \in H , k,k' \in K   \hskip 0.1cm   \hskip 0.1cm  kh = h'k' \}$
 
 Following  N. Andruskiewitsch  and S.Natale's work \cite{AN2}, we are able to define two  double groupoid structures   $\mathcal T$ and $\mathcal T'$. Let's define   "horizontal" and "vertical" products on the squares of $\mathcal T$: 
 \vskip 0.2cm
 let  $ \eenmatrixbase{a}{b}{c }{d}$ and $ \eenmatrixbase{a'}{b'}{c' }{d'}$ be  in $\mathcal T$, they will be composable for the horizontal product $\underset K {\overset h \star}$ if and only if $b=c'$ and $ \eenmatrixbase{a}{b}{c }{d}  \underset K {\overset h \star}\eenmatrixbase{a'}{b'}{b }{d'} = \eenmatrixbase{aa'}{b'}{c }{dd'}$, they will be composable for 
 \vskip 0.1cm
 the vertical product $\underset H {\overset v \star}$ if and only if $d=a'$ and 
 
 $ \eenmatrixbase{a}{b}{c }{d}  \underset H {\overset v \star}\eenmatrixbase{d}{b'}{c' }{d'} = $ \hskip 0.6cm$ \eenmatrixbis{a}{bb'}{cc'}{d'}$. 
 
 One easily sees that $\mathcal T$ is a groupoid  with basis  $H$ for the horizontal product and a groupoid with basis $K$ for the vertical one, of course there  is analogue  properties  for $\mathcal T'$, but its structures also comes from the transpose map: $ \mathcal T \mapsto \mathcal T'$ defined by: $$ \eenmatrixbase{a}{b}{c }{d}  \mapsto   \Big( \eenmatrixbase{a}{b}{c }{d} \Big )^t =  \eenmatrixbase{ c}{d}{a }{b } $$.
 
 \subsubsection{\bf{Remark}}
\label{noel}
Due to lemma \ref{muet}, the { \bf corner maps} associated with these double groupoids (paragraph 1.4 of \cite{AN2}) are all constant and equal to $|S|$.

\subsubsection{\bf{Notation}}
\label{inverse}
For any $t \in \mathcal T$ (resp. $t' \in \mathcal T'$) let's denote by $t^{-h}$ (resp. $t'^{-h}$) and $t^{-v}$ (resp. $t'^{-v}$) the inverse of $t$ (resp.$t'$) for the horizontal and vertical product  respectively. Let's also denote $t^{-hv}$ the double inverse that is $t^{hv} = ({t^{-h}})^{-v} = ({t^{-v}})^{-h}$ and let's denote $t'^{-hv}$ the similar object for   $t'$.

Let $\mathbb C\mathcal T$ (resp. $\mathbb C\mathcal T'$)    be the $\mathbb C$ vector space with canonical basis $\mathcal T$ (resp. $\mathcal T'$)  then  the horizontal product on $\mathcal T$ gives a natural structure  of $*$-algebra to $\mathbb C\mathcal T$  (resp. $\mathbb C\mathcal T'$) with a canonical duality bracket given for any $x \in \mathcal T$ and $x' \in \mathcal T'$ by:
\vskip 0.3cm
 $$<x, x'> = \   \ \left \{ \begin{array}{rl}
&|S|  \  \  \mathrm{if}  \  \  \  x' = x^t \hskip 0.5cm   \\
&  0 \  \ \mathrm{otherwise} \hskip 0.5cm  
\end{array}
\right. 
$$
\vskip 0.9cm

\subsubsection{\bf{Lemma}}
\label{lafinarrivebis}
{\it For any $h,h',h_1,h'_1 \in H$ and $k,k',k_1,k'_1 \in K$, one has:
$$ \eenmatrixbase{k_1}{h_1}{h'_1 }{k'_1}  = (\eenmatrixbase{h}{k}{k' }{h'})^t  \ \ \  \mathrm{{\it if\ \ and \ \ only\ \ if}}  \left \{ \begin{array}{rl}
&hk = k_1h_1\hskip 0.5cm   \\
& h = h'_1  \hskip 0.5cm \\
& k' = k _1 \hskip 0.5cm 
\end{array}
\right.$$ }
\vskip 0.3cm
\begin{dm}
If $ \eenmatrixbase{k_1}{h_1}{h'_1 }{k'_1}  = \eenmatrixbase{k'}{h'}{h}{k} $ \ \ \  \ then $ \left \{ \begin{array}{rl}
&hk = h'_1k'_1= k_1h_1\hskip 0.5cm   \\
& h = h'_1  \hskip 0.5cm \\
& k' = k _1 \hskip 0.5cm 
\end{array}
\right.$
\vskip 0.6cm
Conversely  if $ \left \{ \begin{array}{rl}
&hk = k_1h_1\hskip 0.5cm   \\
& h = h'_1  \hskip 0.5cm \\
& k' = k _1 \hskip 0.5cm 
\end{array}
\right.$ then $ \left \{ \begin{array}{rl}
&hk = k_1h_1 = h'_1k'_1\hskip 0.5cm   \\
&k'h' = hk=  k_1h_1\hskip 0.5cm   \\
& h = h'_1  \hskip 0.5cm \\
& k' = k _1 \hskip 0.5cm 
\end{array}
\right.$
\vskip 0.6cm
so $ \left \{ \begin{array}{rl}
&k = k'_1\hskip 0.5cm   \\
&h' = h_1\hskip 0.5cm   \\
& h = h'_1  \hskip 0.5cm \\
& k' = k _1 \hskip 0.5cm 
\end{array}
\right.$ that means $ \eenmatrixbase{k_1}{h_1}{h'_1 }{k'_1}  = \eenmatrixbase{k'}{h'}{h}{k} $ \ \ .
\end{dm}

\vskip 2cm

\subsubsection{\bf{Theorem (see also \cite{AN2})}}
\label{lafinarrive}
{\it The bracket below   gives to  $\mathbb C\mathcal T$  and $\mathbb C\mathcal T'$  structures of $*$- quantum groupoids in duality, for any $t \in \mathcal T$, one has:
\vskip 0,1cm
\begin{align*}
\Gamma(t)  
&=  \frac{1}{| H \cap K |} \underset{ t_2\underset H{\overset v\star} t_1 = t} \sum t_1  \otimes  t_2 
\\
\kappa( t) &=   t^{-hv}
\\
 \epsilon(t) &=  \   \ \left \{ \begin{array}{rl}
& |H \cap K|  \  \  \mathrm{if}  \  \  \  t  \  \   \mathrm{is\  \ of \ \ the \ \ form} \  \  \eenmatrixbase{h}{e}{e }{h}\hskip 0.5cm   \\
&  0 \  \ \mathrm{otherwise} \hskip 0.5cm  
\end{array}
\right.  
\end{align*}}
\vskip 0.5cm
\begin{dm}  This is simple calculations.
\end{dm}

\subsubsection{\bf Remarks}
\label{inconnue}
Due to remark \ref{noel}, these structures  are  in tight relation with  \cite{AN2} 2.1. The bracket between  $\mathbb C \mathcal T$ and $\mathbb C \mathcal T'$ allows us to define a canonical left action of  the quantum groupoid $\mathbb C \mathcal T$ on the von Neumann algebra $\mathbb C \mathcal T'$.

For any $x  \in \mathcal T$ and any $x'_1,x'_2 \in \mathcal T'$, let's define $x'_1 \triangleright x$  to be the unique element in $\mathcal T$  such that :
$$< \Gamma(x), x'_2 \otimes x'_1> =  <x,x'_2\underset H {\overset h \star}x'_1> =  < x'_1 \triangleright x, ,x'_2>$$

It's very easy to see that for any $y \in \mathcal T'$ and $x \in \mathcal T$:
$$
y \triangleright x =  \   \ \left \{ \begin{array}{rl}
& (y^t)^{-v}\underset H {\overset v \star}x \  \  \mathrm{if  \  \  this \  \ product \  \  exists}   \\
&  0 \  \ \mathrm{otherwise} \hskip 0.5cm  
\end{array}
\right.  
$$

We can give an operator algebra interpretation of these structures. 
\vskip 1.2cm 

\subsubsection{\bf{Proposition}}
\label{calmos}
{\it The  $*$-algebras  $ (C(K) \underset {\rho} \rtimes  S)\underset  {\sigma^I} \rtimes H$  and $  (C(H) \underset {\rho} \rtimes S)\underset  {\sigma^J} \rtimes K $ are respectively isomorphic to the $*$-algebras $\mathbb C \mathcal T$ and $\mathbb C \mathcal T'$.}
\newline

\begin{dm} For any $(h,k,s ) \in H \times K \times S$, let $t(V_h\chi_k\rho(s)) =  \  \    h \underset I\triangleright' k   \eenmatrix{h}{ks}{ }{(h \underset I \triangleleft'  k)s  } \   \ $, for any $(h',k',s' ) \in H \times K \times S$, one has : 

$V_h\chi_k\rho(s)V_{h'}\chi_{k'}\rho(s') = \delta_{ks,h' \underset I  \triangleright' k'}V_{hh'}\chi_{k's^{-1}}\rho(ss')$, and $(V_h\chi_k\rho(s))^* = \rho(s^{-1})\chi_kV_{h^{-1}} = V_{h^{-1}}\chi_{(h  \underset I \triangleright' ks)}\rho(s^{-1})$, from this we can deduce that:

 $t((V_h\chi_k\rho(s))(V_{h'}\chi_{k'}\rho(s'))) = t(V_h\chi_k\rho(s)) \underset H{\overset h \star}t(V_{h'}\chi_{k'}\rho(s'))$ and also:
 \vskip 0.3cm
\begin{align*}
t((V_h\chi_k\rho(s))^*)&= t(\rho(s^{-1})\chi_{h \underset I \triangleleft'  ks}V_{h^{-1}})  
=\hskip 1cm\eenmatrix{h^{-1}}{(h  \underset I \triangleright' ks)s^{-1}}{ks }{ (h \underset I \triangleleft'  ks)^{-1}  }\\
&=   \hskip 1.5cm \eenmatrix{h^{-1}}{(h  \underset I \triangleright' k )}{ks }{ (h \underset I \triangleleft'  ks)^{-1}  }
 \\
&= t((V_h\chi_k\rho(s))^{-h}
\end{align*}
The proposition follows
\end{dm} 
\vskip 1.2cm
We shall see in next chapter, using  a suitable inclusion of von Neumann algebras,  that   the von Neumann algebra crossed product $\mathbb C \mathcal T \rtimes \mathbb C \mathcal T'$ is isomorphic to $\mathbb C[H \cap K] \otimes \mathcal L(\mathbb C^{|H||K|})$.

 \subsection{ ( \cite{BH} chap 4 or \cite{HS}) Quantum groupoids associated with inclusions of von Neumann algebras  coming from relative  matched pairs $H,K$ }
 \label{courage}
 \vskip 0.6cm

 In \cite{BH} is given a very deep study of inclusions of the form $R^H \subset R \rtimes K$ where $H$ and $K$ are {\bf any} subgroups of a group $G$ acting properly and outerly on the  hyperfinite type $II_1$ factor $R$, in such a way that we can identify $G$ with a subgroup of $Out R$,  in particular there is here no ambiguity for the inclusion $H \cap K  \subset Out R$: let's call $\alpha$ the action of $K$ and $\beta$ the one of $H$  here these actions coincide on $S = H \cap K$.    In \cite{BH},   it is proved  that this inclusion is finite depth if and only if the group generated by $H$ and $K$ in $Out R$ is finite,  and, in that situation,  it is irreducible and depth two when $H,K$ is a  matched pair. In this section, using any relative   matched pair, we obtain still depth two but no more irreducible inclusions, so  using \cite{NV2} or \cite{EV} or \cite{D}, these inclusions come from quantum groupoids actions. 
 
   Let's give some facts about these inclusions of the form $R^H \subset R \rtimes K$. First, one can observe that  $Out R$, can be identified, using Sauvageot Connes fusion multiplication and the contragredient procedure,  to a group of $R-R$ bimodules over $L^2(R)$;  and using the sum operation on bimodules,  $Out R$ has also a second  operation with a distributivity property. We shall follow the same notations as in \cite{BH}  and we identify  any element of $Out R$ to the bimodule associated with this element. Let $\gamma =  _{R^H}L^2(R)_{R}$ and $\chi = _RL^2(R \rtimes K)_{R \rtimes K}$, then for any $h \in H$ and $k \in K$ one has:
 
 \begin{align}
 \gamma h = \gamma \  \  \  \  \  \     \  \  k \chi = \chi  \\
 \overline \gamma\gamma = \underset {h \in H} \oplus h \  \  \  \  \  \  \  \   \    \chi\overline \chi  = \underset {k \in K} \oplus k
 \end{align}

\subsubsection{\bf{Lemma}}
\label{jaifaim}
{\it    Let $H,K$ be two finite subgroups of  $Out R$ such that $HK = \{hk/ h\in H , k\in K\}$ is a subgroup of $Out R$, then for any $g$ in $G =HK$, one has $\gamma g\chi = \gamma\chi$.}
\begin{dm}
For any $g$ in $G =HK$, there exist $h \in H$ and $k \in K$ such that $g =hk$, hence due to (1), one has $\gamma g\chi = \gamma hk\chi = \gamma \chi$.
\end{dm}

\subsubsection{\bf{Lemma}}
\label{caen}
{\it  Let $M_0 \subset M_1 \subset M_2 \subset...$ a Jones tower of type $II_1$ factors  such that  the  bimodule $_{M_0}L^2(M_{2})_{M_1}$  is an ampliation of  $_{M_0}L^2(M_{1})_{M_1}$, i.e there exists an integer $n$ such that $_{M_0}L^2(M_{2})_{M_1}$ is isomorphic to $_{M_0\otimes 1}(L^2(M_{1})   \otimes \mathbb C^{n})_{M_1\otimes 1}$ then the inclusion $M_0 \subset M_1$ is depth two.}
\vskip 0.5cm
\begin{dm} In this lemma's conditions, let $x$ be any irreducible $M'_0$-$ M_1$ subbimodule of $\rho = _{M_0}L^2(M_{2})_{M_1}$: it appears $n$ times 
as an irreducible  subbimodule of  \\
$\rho \bar \rho \rho = _{M_0}L^2(M_{1})_{M_1}$, let $x_1,x_2,...x_n$ these subbimodules.  
In   Bratelli's diagramm, there exists an irreducible  subbimodule  $y$ of $\rho \bar \rho = _{M_0}L^2(M_{1})_{M_0}$   which has to be deleted in the principal graph and is connected to $x$. Due to   the Frobenius reciprocity theorem, $y$ is connected to $x_1,x_2,...x_n$,  hence any of these  has to be deleted in the principal graph; but any irreducible subbimodule of $_{M_0}L^2(M_{1})_{M_1}$ is of this form for a good choice of $x$. So the principal graph stops at level two, which is the definition of depth two.\end{dm}

\subsubsection{\bf{Theorem}}
\label{photos}
{\it Let $H,K$ be two finite subgroups of  $Out R$ such that $HK = \{hk/ h\in H , k\in K\}$ is a subgroup of $Out R$, then the inclusion $R^H \subset R \rtimes K$ is depth two. Let $M_2$ be the third element of Jones's tower of inclusion $R^H \subset R \rtimes K$, then there exists a quantum groupoid structure on $ (R^H)' \cap M_2$ over the basis $(R^H)' \cap R \rtimes K$ and an action $\gamma$ of $ (R^H)' \cap M_2$ on $R \rtimes K$ in such a way that the inclusion $R^H \subset R \rtimes K \subset M_2$ is isomorphic to $(R \rtimes K)^{\gamma} \subset R \rtimes K \subset (R \rtimes K)\underset \gamma \rtimes ((R^H)' \cap M_2)$.}
\vskip 0.5cm
\begin{dm}
Let $(M_k)_{k \in K}$ be   Jones tower of   inclusion $R^H \subset R \rtimes K$, due to chap.2 and 3 of \cite{BH}, the  bimodule $_{M_0}L^2(M_{2})_{M_1}$  is equal to $\gamma(\chi\overline \chi\overline \gamma \gamma) \chi$, but using (4) and lemma \ref{jaifaim}, one has

\begin{align*}
\gamma(\chi\overline \chi\overline \gamma \gamma) \chi
&=   \gamma ( \underset{h \in H,k\in K} \oplus hk) \chi = |S| \gamma ( \underset{g \in HK} \oplus g) \chi \\
&= |S| |HK| \gamma \chi = |H||K|\gamma \chi
\end{align*}

So, due to lemma \ref{caen}, the inclusion $R^H \subset R \rtimes K$ is depth two and one can apply  theorem \ref{tenir} or \cite{D} chap 3 or \cite{E1} theorem 9.2 to conclude.
\end{dm}

\vskip 1cm
\subsubsection{\bf{Corollary}}
\label{indiana}
{\it   In the conditions of theorem \ref{photos} the von Neumann algebra  $M'_0 \cap M_3$ is isomorphic to $\mathcal L(H \cap K) \otimes \mathcal L( \mathbb C^{|H||K|})$}.
\label{valentin}
\vskip 0.5cm
\begin{dm}
We have seen that $\gamma(\chi\overline \chi\overline \gamma \gamma) \chi =
|H||K|\gamma \chi
$, which proves the corollary.
\end{dm}

\section{The C*-quantum groupoid structure associated with inclusions of the form $R^H \subset R \rtimes K$}

Till the end of this section, we deal with a finite relative  matched pair $H,K$ acting properly and outerly on $R$. In \cite{HS}, a complete description of Jones tower for the inclusion $R^H \subset R \rtimes K$ is given in the particular case when $H,K$ is a finite  matched pair of groups, a good part of this work remains true for a relative  matched pair, but one must keep in mind that the actions considered in \cite{HS}, are not exactly the one we use but are tightly related: the action of $H$ on $K$ is for example given by $h {\bf .} k = (h \triangleright k^{-1})^{-1}$. To be short, we shall denote $M_0 = R^H$, $M_1 = R \rtimes K$, and $M_2$ will be  the  third element of the basic construction $M_0 \subset M_1 \subset M_2$. All the crossed product here will be defined  "in the right manner", for example $R \rtimes K$ is  the vector space (which is a $*$-algebra)  in $R \otimes \mathcal L(l^2(K))$ generated by the products  $\alpha(r)(1 \otimes \rho(k))$, where $\alpha(r)$ is $r$ viewed in the crossed product as the fonction:  $ k' \mapsto \alpha_{k'}(r)$ i.e $\alpha(r) = ( \underset {k' \in K} \sum \alpha_{k'}(r) \otimes \chi_{k'})$, where $\chi_{k'}$ is the characteristic fonction of $\{k' \}$.

\subsection{ The $^*$-algebra structure of $M'_0 \cap M_2$ }
\label{temoin}
One can easily see what are the basis $(R^H)' \cap R \rtimes K$ of quantum groupoid $M'_0 \cap M_2$.
\subsubsection{\bf{Lemma(BH or HS)} }
\label{dietmar}
{\it  The algebra $(R^H)' \cap R \rtimes K$ is isomorphic to 
the group algebra $\mathcal  L (S)$. }
\newline
\begin{dm}
If $u_k$ and $v_h$ are canonical implementations of $\alpha$ and $\beta$ on $L^2(R)$, one can suppose $u_x = v_x$ for any $x$ in $S$, these $u_x$ generate a *-algebra isomorphic to $\mathbb C[S]$ and are clearly in $(R^H)' \cap R \rtimes K$, on the  other hand, using the computation in the proof of 4.1 in \cite{BH}, one has: $dim ((R^H)' \cap R \rtimes K )= card S$. The lemma follows.
\end{dm}

\subsubsection{\bf{Remark}}
\label{photocopie}
The basis $(R^H)' \cap R \rtimes K$ of quantum groupoid $(R^H)' \cap R \rtimes K$ are commutative if and only if $S$ is abelian. Hence, when $S$ is non abelian,  this  quantum groupoid structure  does not come from a  matched pair of groupoids ( see \cite{Val2}). 
\vskip 1cm
Now let's give a description of the two  first steps of the basic construction for the inclusion $R^H \subset R \rtimes K$. For our computations, it will be more convenient to   express as soon as possible  {\bf all the algebras in $M_3$}. As there is two actions, namely $\alpha$ and $\beta$ of respectively $K$ and $H$ on $R$, there is also   two actions $\alpha^1$ and $\beta^1$ of $K$ and $H$ on $R^1 = R \otimes \mathcal L(l^2(K))$. As well known (see \cite{Val3} th 5.3 for a very old reference...), one can identify $R^1$ with the double crossed product $R \underset \alpha \rtimes K\underset {\hat \alpha} \rtimes \hat K$,   $\alpha^1$ is the bidual action $  \hat {\hat \alpha}$,  so $\alpha^1 = \alpha \otimes Ad \lambda$ the fixed points  algebra of which is $M^1$, 
and   $\beta^1$ is the action of $H$ on  $R^1 = R \otimes \mathcal L(l^2(K))$ defined by $\beta^1 = \beta \otimes 1 $.  Let $(w_h)_{h \in H}$ be a group of unitaries of $R^1\underset{\beta^1} \rtimes H$ implementing $\beta^1$, hence  $R^1$ is the set of  sums : $\underset{ k,k' \in K} \sum x_{k,k'}  \otimes \rho(k)\chi_{k'}$, where $x_{k,k'}$ is in $R$, and $R^1\underset{\beta^1} \rtimes H$ (viewed in $L^2(R^1)$) is the set of sums: $\underset{h \in H; k,k' \in K} \sum (  x_{k,k',h} \otimes \rho(k)\chi_{k'})w_h$ , where $x_{k,k'}$ is a non zero element in $R$, $\chi_k$ is the multiplication operator by the characteristic fonction of $\{k\}$. In fact, for any algebraic basis $(x_i)_{i \in I}$ of $R$ the family $\big (( x_i \otimes \rho(k)\chi_{k'})w_h\big )_{i \in I, h  \in H; k,k' \in K} $ is an algebraic basis for $R^1\underset{\beta_1} \rtimes H$, also  the family $\big (( x_i \otimes \lambda(k)\chi_{k'})w_h\big )_{i \in I, h,  \in H; k,k' \in K} $ is an algebraic basis for $R^1\underset{\beta_1} \rtimes H$. 
Let $\tau$ be the normalized tracial state for the Jones 's tower, and $\tilde \tau $ the one of $R^1\underset{\beta_1} \rtimes H$, let $E: R^1\underset{\beta_1} \rtimes H \to R  \rtimes K$ be the $\tilde \tau$ preserving conditional expectation, then by routine calculations, for any $k,k' \in K$, $x \in R$ and $h \in H$, one has: 

$\tilde \tau (\alpha(x)(1 \otimes \rho(k)\chi_{k'})w_h) = \frac{1}{|K|}\tau(x)\delta_{k,e}\delta_{h,e} $ 

 $E(\alpha(x)(1 \otimes \rho(k)\chi_{k'})w_h) = \frac{1}{|K|} \alpha(x)(1 \otimes \rho(k))\delta_{h,e},$ where $\alpha(x)$ is $x$ viewed in $M_1 = R  \rtimes K$, i.e $\alpha(x) =  \underset{k_0} \sum \alpha_{k_0}(x) \otimes \chi_{k_0} .$

\vskip 1cm

\subsubsection{\bf{Proposition(HS)}}
\label{corrompu}
{\it With Jones's notations, $(M_2,e_1)$  can be identified with $(R\otimes  \mathcal L(l^2(K))\underset{\beta^1}\rtimes H, \frac{1}{|H|}(1 \otimes \chi_e)\underset{h \in H} \sum w_h)$).}
\vskip 0.8cm
\begin{dm} Let $e$ be the projection  in $M_0$ ($= R^H$) equal to  $ \frac{1}{|H|}(1 \otimes \chi_e)\underset{h \in H} \sum w_h$. For any $r$   in $M_0$, let $ \alpha(r)$ be $r$  viewed in $R^1\underset{\beta^1} \rtimes H$ (as a von Neumann subalgebra of $R \otimes \mathcal L(l^2(K))$ i.e:  $\alpha(r)= \underset{k_0} \sum \alpha_{k_0}(r) \otimes \chi_{k_0}$. 

As $r$ is a fixed point for $\beta$, one has:
\begin{align*}
&e\alpha(r)  - \alpha(r)e 
= \\
=& \frac{1}{|H|}(1 \otimes \chi_e)\underset{h \in H} \sum w_h\underset{k_0} \sum (\alpha_{k_0}(r) \otimes \chi_{k_0}) -(\underset{k_0} \sum \alpha_{k_0}(r) \otimes \chi_{k_0})\frac{1}{|H|}(1 \otimes \chi_e)\underset{h \in H} \sum w_h \\
=& \frac{1}{|H|}\underset{h \in H} \sum w_h  (r  \otimes \chi_e) - (r  \otimes \chi_e)w_h  
\end{align*}

So we have $[e,\alpha(r)] = 0$, and:

\begin{align*}
E(e) 
&= \frac{1}{|H||K|} 1= [R^1\underset{\beta^1} \rtimes H : R^1]^{-1}[R \underset \alpha \rtimes K\underset {\hat \alpha} \rtimes \hat K: R \underset \alpha \rtimes K]^{-1}1 \\
&= [R^1\underset{\beta_1} \rtimes H : M_1]^{-1}1
\end{align*}
and also:
\begin{align*}
E(e) 
&= \frac{1}{|H||K|} 1= [R :R^H]^{-1}[M_1 : R]^{-1}1 = [M_1:M_0]^{-1}1\\
&=  [M_2:M_1]^{-1}1
\end{align*}

Hence, one can apply the equivalence $1^0$ and $ 2^0$ of Proposition 1.2 in \cite{PP} to conclude that $(M_2,e_1)$ can be identified with $(R\otimes  \mathcal L(l^2(K))\underset{\beta^1}\rtimes H, \frac{1}{|H|}(1 \otimes \chi_e)\underset{h \in H} \sum w_h)$).
\end{dm}
\vskip 1.5cm

Let's give a description of the second step of the basic construction for the inclusion $R^H \subset R \rtimes K$. One can consider the double crossed product $(R^1\underset{\beta^1} \rtimes H )\underset{\hat \beta^1} \rtimes \hat H$ ( $R^1 =  R \otimes  \mathcal L (l^2(H))$) which can be identified with $R^2 = R^1 \otimes \mathcal L (l^2(H)) = R \otimes  \mathcal L (l^2(K)) \otimes  \mathcal L (l^2(H))  $, we can define an  action on $R^2$ by $\alpha^2 = \alpha^1 \otimes 1 = \alpha \otimes Ad \lambda \otimes 1$   and consider a group of unitaries  $(v_k)_{k \in K}$ implementing $\alpha^2$ such that the  crossed product $R^2\underset{\alpha^2} \rtimes K = (R \otimes \mathcal L (l^2(K)  \otimes \mathcal L (l^2(H)))\underset{\alpha^2} \rtimes K$ has a canonical basis of the form $\big (( y_j \otimes \rho(h)\chi_{h'})v_k\big )_{i \in J h, h' \in H; k' \in K} $ for any basis  $(y_j)_{j \in J}$ of $R^1 =  R \otimes  \mathcal L (l^2(H))$.
\vskip 0,8cm
\subsubsection{\bf{Proposition(HS)}}
\label{corrompubis}
{\it With Jones 's notations, $(M_3,e_2)$ can be identified with $(R\otimes  \mathcal L(l^2(K))\otimes  \mathcal L(l^2(H))\underset{\alpha^2}\rtimes K, \frac{1}{|K|}(1 \otimes 1 \otimes \chi_e)\underset{k \in K} \sum v_k)$ .}
\vskip 0.8cm
\begin{dm}  As $M_1$ is just the fixed points algebra of $\alpha^1$, the demonstration is the same as proposition \ref{corrompu}
\end{dm}

\vskip 1,2cm

Let's  see, with our identifications, how $M_2$ is included in $M_3$:   any $z \in M_2$, belonging to $R \otimes  \mathcal L (l^2(H))$, must be viewed in $M_3$ as the operator $\underset {h \in H} \sum \beta^1_h(z) \otimes  \chi_h$ in $R \otimes \mathcal L(K) \otimes \mathcal L(H)$, so any element of $M_2$ of the form $1 \otimes x$ where $x$ is in $\mathcal L (l^2(K))$,  is seen as $1 \otimes x \otimes 1$ in $M_3$  and $ w_s$ is just the operator $1 \otimes 1 \otimes \rho(s)$. One can identify  $\tilde \tau$ with the "Markov trace" $\tau$. If $M_3$ is viewed as a triple crossed product $((R^1\underset{\beta^1} \rtimes H) \underset{\hat \beta^1} \rtimes \hat H )\underset{\alpha^2} \rtimes K$, for any element $r^1 \in R^1$ $h,h' \in H$ and $k \in K$: $\tau(\beta^1(r^1)(1 \otimes 1 \otimes \rho(h)\chi_{h'})v_k) =  \frac{1}{|H|}\tau(r^1)\delta_{h,e}\delta_{k,e}. $.

\subsubsection{\bf{Notations}}
\label{fievre}
For any $ h \in H,  k \in K,  k' \in p_2(k^{-1}h^{-1}k)$,  let's define $w_{k,h,k'} =   1 \otimes \rho(k')\chi_{k}\otimes \rho(h)$ .

\subsubsection{\bf{Lemma}}
\label{graphe}
{\it i)   $ \{ w_{k,h,k'}  /  h \in H,  k \in K,  k' \in p_2(k^{-1}h^{-1}k) \} $ defines a basis of $ M'_0 \cap M_2 $.

ii) for any $ h \in H,  k \in K,  k' \in p_2(k^{-1}h^{-1}k)$, $ w_{k,h,k'}$ is a partial isometry with initial  support $ 1 \otimes \chi_k \otimes 1$  and final support $1 \otimes \chi_{kk'^{-1}} \otimes 1$.}
\vskip 0.5cm
\begin{dm}
Let's make the computations in $M_2$. Let $r$ be any element in $M_0 = R^H$ and $ y = \underset{k_0} \sum \alpha_{k_0}(r) \otimes \chi_{k_0} $ the same viewed in $M_2 = R^1\underset{\beta_1} \rtimes H$, for any element $x = \underset{h \in H; k,k' \in K} \sum (  x_{k,k',h} \otimes \rho(k)\chi_{k'})w_h$, x commutes with $M_0$ means that for any $r \in M_0$, one has: 

\begin{multline*} \big (\underset{k_0} \sum \alpha_{k_0}(r) \otimes \chi_{k_0}\big )\underset{h \in H; k,k' \in K} \sum (  x_{k,k',h} \otimes \rho(k)\chi_{k'})w_h  = \\
\underset{h \in H; k,k' \in K} \sum (  x_{k,k',h} \otimes \rho(k)\chi_{k'})w_h \big (\underset{k_0} \sum \alpha_{k_0}(r) \otimes \chi_{k_0}\big )
\end{multline*}

On the one hand, one has:
\begin{align*}
 \big (\underset{k_0} \sum \alpha_{k_0}(r) \otimes \chi_{k_0}\big )&\underset{h \in H; k,k' \in K} \sum (  x_{k,k',h} \otimes \rho(k)\chi_{k'})w_h 
 =  \\
&=  \underset{h \in H; k_0,k,k' \in K} \sum (  \alpha_{k_0}(r)x_{k,k',h} \otimes \chi_{k_0} \rho(k)\chi_{k'})w_h \\
&=  \underset{h \in H; k_0,k,k' \in K} \sum (  \alpha_{k_0}(r)x_{k,k',h} \otimes \rho(k)\chi_{k_0k} \chi_{k'})w_h \\
&=  \underset{h \in H; k_0,k \in K} \sum (  \alpha_{k_0}(r)x_{k,k_0k,h} \otimes \rho(k)\chi_{k_0k} )w_h \\
&=  \underset{h \in H; k_0,k \in K} \sum (  \alpha_{k_0k^{-1}}(r)x_{k,k_0,h} \otimes \rho(k)\chi_{k_0} )w_h 
\end{align*}

 On the other  hand, one has:
\begin{align*}
\underset{h \in H; k,k' \in K} \sum (  x_{k,k',h} \otimes \rho(k)\chi_{k'})w_h &\big (\underset{k_0} \sum \alpha_{k_0}(r) \otimes \chi_{k_0}\big )
 =  \\
&= \underset{h \in H; k,k' \in K} \sum (  x_{k,k',h} \otimes \rho(k)\chi_{k'})\beta^1_h (\underset{k_0} \sum \alpha_{k_0}(r) \otimes \chi_{k_0})w_h 
\\
&= \underset{h \in H; k,k' \in K} \sum (  x_{k,k',h} \otimes \rho(k)\chi_{k'})(\underset{k_0} \sum \beta_h (\alpha_{k_0}(r) )\otimes \chi_{k_0})w_h 
\\
&=   \underset{h \in H; k_0,k,k' \in K} \sum (  x_{k,k',h} \beta_h (\alpha_{k_0}(r) )\otimes \rho(k)\chi_{k'})\chi_{k_0})w_h \\
&=   \underset{h \in H; k_0,k\in K} \sum (  x_{k,k_0,h} \beta_h (\alpha_{k_0}(r) )\otimes \rho(k)\chi_{k_0})w_h 
 \end{align*}

So $x$ is in $M'_0 \cap M_2$ if and only if for any $h \in H; k_0,k\in K$ and any $r \in R^H$, one has:
$  \alpha_{k_0k^{-1}}(r)x_{k,k_0,h}  =  x_{k,k_0,h} \beta_h (\alpha_{k_0}(r) )$. Applying $ \alpha_{k{k_0}^{-1}}$, this is equivalent to:
$$ r \alpha_{k{k_0}^{-1}}(x_{k,k_0,h} ) =    \alpha_{k{k_0}^{-1}} (x_{k,k_0,h})(\alpha_{k{k_0}^{-1}} \beta_h \alpha_{k_0})(r) $$

So, using lemma 3.1 of \cite{HS}, for any $h \in H; k_0,k\in K$  such that $ x_{k,k_0,h} \not = 0$, there exist $\lambda$ in $\mathbb C - \{0\}$, $v_{h,k_0,k}$ in $\mathcal U(R)$ (the set of unitaries of $R$)  and $h'$ in $H$ such that:
 $ \alpha_{k{k_0}^{-1}} (x_{k,k_0,h}) = \lambda v_{h,k_0,k}$ and $Ad(v_{h,k_0,k})\alpha_{k{k_0}^{-1}} \beta_h \alpha_{k_0} = \beta_{h'}$. Hence one has: $Ad(v_{h,k_0,k}) = \beta_{h'}\alpha_{{k_0}^{-1}} \beta_{h^{-1}} \alpha_{k_0}\alpha_{k^{-1}}$, this equality in $G$ means: $v_{h,k_0,k} \in Z(R) (= \mathbb C)$ and $h'{k_0}^{-1}h^{-1}k_0k^{-1} = e$, which are exactly the two conditions:  $x_{h,k_0,k} \in  \mathbb C$ and $k \in p_2(k_0^{-1}h^{-1}k_0)$.

ii) This is an easy computation.

\end{dm}

\vskip 0.8cm

\subsubsection{\bf{Corollary }}
\label{complement}
{\it i) The    family $(1 \otimes \rho(s)\chi_k \otimes 1)_{s \in H \cap K , k \in K}$ is a linear basis for a $*$-subalgebra of $M'_0 \cap M_2$ isomorphic to  the crossed product $  C(K) \underset {\rho} \rtimes (H \cap K)$, 

ii) the family $ (1 \otimes \rho(s)\otimes 1)_{s \in  H \cap K}$ is a group of unitaries  and a linear basis for $M'_0 \cap M_1$,

iii) the family $ (1 \otimes \lambda(s)\otimes \rho(s))_{s \in H \cap K }$ is a group of unitaries  and a linear basis for $M'_1 \cap M_2 $.}
\vskip 0.2cm
\begin{dm} The assertion i) is trivial. 
For any $s \in S$ and $k \in K$, one has: $1 \otimes \rho(s) =  \underset {k_1 \in K} \sum w_{k_1,e,s}$ and $1 \otimes \chi_k = w_{k,e,e}$.
For any $s \in S$, one has: $ \{ 1 \otimes \rho(s) / s \in S\} \subset M_1 \cap (M'_0 \cap M_2) = M'_0 \cap M_1$ and for dimension reasons ii) follows. As the families $ (1 \otimes \lambda(s))_{s \in S }$ and $ (w_s)_{s \in S }$ are commuting groups of unitaries, so $ ((1 \otimes \lambda(s))w_s)_{s \in S }$ is a group of unitaries, and one has:  $\lambda(s) =  \underset {k \in K} \sum \rho(k^{-1}s^{-1}k)\chi_k$, hence:  $(1 \otimes \lambda(s))w_s = \underset {k \in K} \sum  w_{k,s,k^{-1}s^{-1}k}$ hence it is in $M_2$. Let  $r$ be any element in $ R$  and $ y = \underset{k_0} \sum \alpha_{k_0}(r) \otimes \chi_{k_0} $ the same viewed in $M_1 = R  \rtimes H$, for any $\sigma \in S$ one has: 
\begin{align*}
y(1 \otimes \rho(\sigma)) w_s(1 \otimes \lambda(s))
&= \underset{k_0} \sum( \alpha_{k_0}(r) \otimes \chi_{k_0} \rho(\sigma))w_s(1 \otimes \lambda(s))  \\
&= w_s(1 \otimes \lambda(s))(\underset{k_0} \sum (\beta_{s^{-1}}\alpha_{k_0})(r) \otimes \chi_{s^{-1}k_0} \rho(\sigma) )\\
&= w_s(1 \otimes \lambda(s))(\underset{k_0} \sum( (\alpha_{s^{-1}}\alpha_{k_0})(r) \otimes \chi_{s^{-1}k_0} \rho(\sigma))  \\
&= w_s(1 \otimes \lambda(s))(\underset{k_0} \sum( (\alpha_{s^{-1}k_0})(r) \otimes \chi_{s^{-1}k_0} \rho(\sigma))\\
&= w_s(1 \otimes \lambda(s))(\underset{k'} \sum( (\alpha_{k'})(r) \otimes \chi_{k'} \rho(\sigma))\\
&= w_s(1 \otimes \lambda(s))y( 1 \otimes  \rho(\sigma))
\end{align*}

So, the unitary $w_s(1 \otimes \lambda(s )$ is in $M'_1 \cap M_2$ and for dimension reasons  iii) follows. 
\end{dm}

\vskip 1.2cm

\subsubsection{\bf{Lemma}}
\label{poubelles}
{\it  Using the notations of paragraph \ref{simplifier}, for any  $h \in H$ and $k \in K$, let $k'_I(k,h)= (h \underset I \triangleright k)^{-1}k$, then $k'_I(k,h) $ is in $p_2(k^{-1}h^{-1}k)$ and for any $ h' \in H$ and $k' \in K$, one has: 
$$k'_I(h' \underset I \triangleright' k',h) k'_I(k',h' ) = k'_I(k',hh')$$}
\vskip 0.3cm
\begin{dm}
For any  $h \in H$ and $k \in K$, one has: $ k^{-1}h^{-1}k = (hk)^{-1}k = (h \underset I \triangleleft k)^{-1}(h \underset I \triangleright k)^{-1}k = (h \underset I \triangleleft k)^{-1}k'_I(k,h)$, hence, $k'_I(k,h)$ is in 
$p_2(k^{-1}h^{-1}k)$.
For any $h,h' \in H$ and $k' \in K$, one has: 
\begin{align*}
k'(k'_I(h' \underset I \triangleright' k',h) k'_I(k',h' ))^{-1}
&=k'(k'_I(k',h' ))^{-1}(k'_I(h' \underset I \triangleright' k',h))^{-1}  \\
&= (h' \underset I \triangleright' k')(k'_I(h' \underset I \triangleright' k',h))^{-1}  = h \underset I \triangleright' (h' \underset I \triangleright' k') \\
&= (hh' )\underset I \triangleright' k'
\end{align*}
So $k'_I(h' \underset I \triangleright' k',h) k'_I(k',h' ) = k'_I(k',hh')$.
\end{dm}

\vskip 1cm

\subsubsection{\bf {Corollary and notations}}
\label{domoterm}
For any  $h \in H$ and $k \in K$, the element 
$    W^I_h = \underset { k \in K}\sum  1 \otimes \rho((h \underset I \triangleright k)^{-1}k)\chi_{k}\otimes \rho(h) $  is in  $M'_0 \cap M_2$.  
\vskip0.2cm
\begin{dm}
This comes from lemmas \ref{graphe} and \ref{poubelles}
\end{dm}
\vskip 0.6cm

\subsubsection{\bf{Theorem}}
\label{banlieue}
{\it The family  $( W^I_h)_{h \in H}$ is a one parameter group of unitaries in $M'_0 \cap M_2$, and   it  implements an action of $H$ on the  $*$-subalgebra of $M'_0 \cap M_2$ generated by the family $(1 \otimes \rho(s)\chi_k \otimes 1)_{s \in S, k \in K}$ which is equivalent to the action $\sigma^I$ defined in theorem \ref{roy}. The $*$-algebra   $M'_0 \cap M_2$ is isomorphic to the double crossed product $ ( C(K) \underset {\rho} \rtimes (H \cap K)) \underset {\sigma^I}\rtimes H $ and to $\mathbb C \mathcal T$. } 
\vskip 0.2cm
\begin{dm}
For any $k_0 \in K$, and $h\in K$,  due to \ref{graphe} ii), one has:
\begin{align*}
  W^I_h(1 \otimes \chi_{k_0}\otimes 1) { W^I_h}^* 
 &=   \underset {k,k^1\in H } \sum w_{k,h,k'_I(k,h)}  (1 \otimes \chi_{k_0})  w^*_{k^1,h,k'_I(k^1,h)} 
\\
 &=    w_{k_0,h,k'_J(k_0,h)}   w^*_{k_0,h,{k'_I}(k_0,h)}  \\
&=  1   \otimes \chi_{ h \underset I\triangleright'  k_0} \otimes 1
\end{align*}
Suming this equality for all $k_0 \in K$ gives that $W^I_{k_0}$ is a unitary.

Now for any $s \in S$, $k \in K$ and $h \in H$, one easily sees that  $k'_I(ks^{-1},h)= sk'_I(h,k)s^{-1}$, from this one  deduces that: 

\begin{align*}
(1 \otimes \rho(s) &\otimes 1)w_{k,h,k'_I(k,h)}
= 1 \otimes\rho(sk'_I(k,h))\chi_{k}\otimes \rho(h) = 1 \otimes  \rho(sk'_I(k,h)s^{-1}s)\chi_{k}\otimes \rho(h)\\
&= 1 \otimes  \rho(h'_I(ks^{-1},h)\rho(s)\chi_{k}\otimes \rho(h)  =  (1 \otimes  \rho(k'_I(ks^{-1},h)\chi_{ks^{-1}} \otimes \rho(h) )( 1 \otimes  \rho(s) \otimes 1)\\
&=  w_{ks^{-1},h,h'_J(ks,h)} ( 1   \otimes \rho(s) \otimes 1)
\end{align*}

From this one deduces that: 
\begin{align*}
 W^I_h(1   \otimes \rho(s)\otimes 1){ W^I_h}^* = 1   \otimes \rho(s)\otimes 1
 \end{align*}

\vskip 0.5cm

Also, for any $h,h'$ in $H$, due to lemma \ref{graphe} ii) and lemma \ref{poubelles}, we have:

\begin{align*}
 W^I_hW^I_{h'}
&=     \underset{k,k' \in  K } \sum      w_{k,h, k'_I(k,h)}w_{k',h', k'_I(k',h')} =     \underset{ k' \in  K } \sum      w_{h' \underset I \triangleright' k',h, k'_I(h' \underset I \triangleright' k',h)}w_{k',h', k'_I(k',h')} \\
&=     \ \underset{ k' \in  K } \sum   (1 \otimes \rho( k'_I(h' \underset I \triangleright' k',h)) \rho( k'_I(k',h' ))\chi_{k'}\rho(h)\rho(h')\\
&=     \ \underset{ k' \in  K } \sum   (1 \otimes \rho( k'_I(h' \underset I \triangleright' k',h) k'_I(k',h' ))\chi_{k'}\rho(hh')
 \\
&=     \ \underset{ k' \in  K } \sum   (1 \otimes \rho( k'_I( k',hh') )\chi_{k'}\rho(hh')= W^I_{hh'}
\end{align*}

So the group of unitaries $( W^I_h)_{h \in I}$ implements an action on the $*$-algebra generated by the families $(1 \otimes \rho(s) \otimes 1)$ and $( 1\otimes \chi_k \otimes 1)$ which is clearly isomorphic to $C(K) \underset {\rho} \rtimes S$ and this action is equivalent to $\sigma^I$.

For any $s  \in S$, $h \in H$ and $k \in K$, one has: $(1 \otimes \rho(s) \otimes 1) W^I_h( 1 \otimes \chi_{k} \otimes 1) = w_{k,h, sk'_j(k,h)}$, this proves that the crossed product  is exactly $M'_0 \cap M_2$ and so, by theorem \ref{roy} , $M'_0 \cap M_2$ is isomorphic to the crossed product $ ( C(K) \underset {\rho} \rtimes S)\underset {\sigma^I_h}\rtimes H $, the theorem follows. \end{dm}
\vskip 1cm

\subsubsection{\bf {Corollary and notations}}
\label{jesouffre}
{\it For any $t \in \mathcal T$, let $(h,k,s)$ be the unique element 
\vskip 0,3cm
in $H \times K\times S$ such that $t = \  \    h \underset I\triangleright' k   \eenmatrix{h}{ks}{ }{(h \underset I \triangleleft'  k)s  } \   \ $, let $\theta_t^I$ be the element in $M'_0 \cap M_2$
\vskip 0,2cm
 equal to $W^I_h(1 \otimes \chi_k\rho(s) \otimes 1)$, then the family $(\theta^I_t)_{t \in \mathcal T}$ is a basis of $M'_0 \cap M_2$ such that for any $t' \in \mathcal T$, one has:

$$\theta_t^I\theta_{t'}^I =  \   \ \left \{ \begin{array}{rl}
& \theta^I_{ t \underset H{\overset h \star}t'} \  \mathrm{if \ \  t \ \ and\ \ t' \ \ are \ \ composable\ \  for} \ \ \underset H{\overset h \star}  \hskip 0.5cm   \\
&  0 \  \ \mathrm{otherwise} \hskip 0.5cm  
\end{array}
\right. $$
$$ ( \theta_t^I)^* = \theta^I_{t^{-h}} $$
}
\newline
\begin{dm}
Using theorem \ref{banlieue}, this is just a reformulation of proposition \ref{calmos}
\end{dm}

\vskip 1.8cm
\subsection{The $*$-algebra structure of  $M'_1 \cap M_3$}
\subsubsection{\bf{Notations}}
\label{fievrebis}
For any $h \in H$, $k \in K$ and $h' \in p'_2(h^{-1}k^{-1}h)$,  let's define $k'= h'h^{-1}kh$  and $v_{h,k,h'} = (1 \otimes \lambda(k'k^{-1}) \otimes \rho(h')\chi_h)v_k$ .

\subsubsection{\bf{Lemma}}
\label{graphebis}
{\it i)   $ \{ v_{h,k,h'}  /  h \in H,  k \in K,  h' \in p'_2(h^{-1}k^{-1}h) \} $ defines a basis of $ M'_1 \cap M_3 $.

ii) for any $ h \in H,  k \in K,  h' \in p'_2(h^{-1}k^{-1}h)$, $ v_{h,k,h'} $ is a partial isometry with initial   support $1 \otimes 1 \otimes \chi_h$ and final support  $    1 \otimes 1 \otimes \chi_{hh'^{-1}}$.}
\vskip 0.5cm
\begin{dm}
Let $y$ be any element in $M_1 $ and  $ \underset{h_1} \sum \beta^1_{h_1}(y) \otimes \chi_{h_1} $ the same viewed in $M_3$. Let $x = \underset{k \in K; h,h' \in H} \sum (  x_{h,h',k} \otimes \rho(h)\chi_{h'})v_k$ be any element of $M_3$, x commutes with $M_1$ means that for any $y \in M_1$, one has: 

\begin{multline*} \big (\underset{h_1} \sum \beta^1_{h_1}(y) \otimes \chi_{h_1}\big )\underset{k \in K; h,h' \in H} \sum (  x_{h,h',k} \otimes \rho(h)\chi_{h'})v_k = \\
\underset{k \in K; h,h' \in H} \sum (  x_{h,h',k} \otimes \rho(h)\chi_{h'})v_k \big (\underset{h_1} \sum \beta^1_{h_1}(y) \otimes \chi_{h_1}\big )
\end{multline*}

On the  one hand, one has:
\begin{align*}
  \big (\underset{h_1} \sum \beta^1_{h_1}(y) \otimes \chi_{h_1}\big )&\underset{k \in K; h,h' \in H} \sum (  x_{h,h',k} \otimes \rho(h)\chi_{h'})v_k = \\
&=  \underset{h_1,h,h' \in H; k \in K} \sum (  \beta^1_{h_1}(y)x_{h,h',k} \otimes \chi_{h_1} \rho(h)\chi_{h'})v_k \\
&=  \underset{h_1,h,h' \in H; k \in K} \sum (  \beta^1_{h_1}(y)x_{h,h',k} \otimes \rho(h)\chi_{h_1h} \chi_{h'})v_k  \\
&=  \underset{h_1,h\in H; k \in K} \sum (  \beta^1_{h_1}(y)x_{h,h_1h,k} \otimes \rho(h)\chi_{h_1h})v_k\\
&=  \underset{h,h_0 \in H; k \in K} \sum (  \beta^1_{h_0h^{-1}}(y)x_{h,h_0,k} \otimes \rho(h)\chi_{h_0})v_k 
\end{align*}

On the other  hand, one has:
\begin{align*}
\underset{k \in K; h,h' \in H} \sum (  x_{h,h',k} \otimes \rho(h)\chi_{h'}&)v_k \big (\underset{h_1} \sum \beta^1_{h_1}(y) \otimes \chi_{h_1}\big )
 =  \\
&=\underset{k \in K; h,h' \in H} \sum (  x_{h,h',k} \otimes \rho(h)\chi_{h'})\alpha^2_k\big (\underset{h_1} \sum \beta^1_{h_1})(y) \otimes \chi_{h_1}\big )v_k \\
&=\underset{k \in K; h,h' \in H} \sum (  x_{h,h',k} \otimes \rho(h)\chi_{h'})\big (\underset{h_1} \sum(\alpha^1_k \beta^1_{h_1})(y) \otimes \chi_{h_1}\big )v_k\\
&=\underset{k \in K; h,h' ,h_1\in H} \sum (  x_{h,h',k} (\alpha^1_k \beta^1_{h_1})(y)\otimes \rho(h)\chi_{h'} \chi_{h_1})v_k\\
&=\underset{k \in K; h,h_1\in H} \sum (  x_{h,h_1,k} (\alpha^1_k \beta^1_{h_1})(y)\otimes \rho(h)\chi_{h_1})v_k 
 \end{align*}

So $x$ is in $M'_1 \cap M_3$ if and only if for any $h ,h_1 \in H; k\in K$ and any $y \in M_1$, one has:
$  \beta^1_{h_1h^{-1}}(y)x_{h,h_1,k}   =   x_{h,h_1,k} (\alpha^1_k \beta^1_{h_1})(y)$
. Applying $ \beta^1_{h{h_1}^{-1}}$, this is equivalent to:
$$ y\beta^1_{h{h_1}^{-1}}(x_{h,h_1,k} ) =    \beta^1_{h{h_1}^{-1}}(x_{h,h_1,k} )(\beta^1_{h{h_1}^{-1}}\alpha^1_k \beta^1_{h_1})(y) $$

As $\alpha^1 = \alpha \otimes Ad\lambda$, $\alpha^1$ is outer and  using lemma 3.1 of \cite{HS}, for any $h,h_1 \in H; k\in K$  such that $x_{h,h_1,k} \not = 0$, there exist $\mu$ in $\mathbb C - \{0\}$, $v_{h,h_1,k}$ in $\mathcal U(R \otimes \mathcal L(l^2(K))$   and $k'$ in $K$ such that:
 $\beta^1_{h{h_1^{-1}}}(x_{h,h_1,k} ) = \mu v_{h,h_1,k}$ and $Ad(v_{h,h_1,k})\beta^1_{h{h_1^{-1}}}\alpha^1_k \beta^1_{h_1} = \alpha^1_{k'}$. Hence one has:  $Ad(v_{h,h_1,k}) = \alpha^1_{k'}\beta^1_{{h_1^{-1}}}\alpha^1_{k^{-1}}\beta^1_{h_1{h}^{-1}}$, which is equivalent to the fact that: for any $r \in R$ and $x \in \mathcal L(l^2(K))$, one has:  
\begin{align*}
Ad(v_{h,h_1,k})(r \otimes 1) 
&= \alpha_{k'}\beta_{{h_1^{-1}}}\alpha_{k^{-1}}\beta_{h_1{h}^{-1}}(r) \otimes 1\\
Ad(v_{h,h_1,k})(1 \otimes x) 
&= 1 \otimes Ad \lambda(k' k^{-1})(x)
 \end{align*}
 this equality  implies  there exists an element $v_{h,h_1,k}^1 \in \mathcal L(l^2(K))$ such that $v_{h,h_1,k} = 1Ê\otimes v_{h,h_1,k}^1 $, $k'h_1^{-1}k^{-1}h_1{h}^{-1}  = e$,  and $Ad(v_{h,h_1,k}^1)=Ad \lambda(k' k^{-1}) $ which are exactly the two conditions:  $x_{h,h_1,k}  \in \mathbb C( 1 \otimes \lambda(k' k^{-1})) $ and $h \in p'_2(h_1^{-1}k^{-1}h_1)$.

ii) This is an easy computation.

\end{dm}

\subsubsection{\bf {Remarks and notations}}
\label{domotermbis}

Using the same argument as  in  lemma \ref{poubelles},    the element ${h'_J}^{-1}(h,k) = (k \underset J \triangleright  h)^{-1}h$ is   in $p'_2(k^{-1}h^{-1}k)$  and  if  $k'$ is the element  in $K$ which enters in the definition of $v_{h,k,h'_J}$ (cf. \ref{fievrebis}), one has  $k' h^{-1}k^{-1}h{h'_J}^{-1}  = e$, so $(k \underset J \triangleright  h)k' = kh$, hence, using notations \ref{andros}, one has:  $k' = k \underset J \triangleleft h$.
 With these notations  one can define an element  $W^J_k$, in  $M'_1 \cap M_3$,  by:

$$   W^J_k = \underset {h\in H } \sum v_{h,k,h'_J(h,k)} = \underset {h\in H } \sum (1 \otimes \lambda((k \underset J \triangleleft h)k^{-1}) \otimes \rho((k \underset J \triangleright  h)^{-1}h))\chi_{h}v_k .$$

\subsubsection{\bf {Lemma}}
\label{bepc}
{\it For any $h \in H$ and $k,k' \in K$, one has: 
$$ h'_J(k \underset J \triangleright  h,k')h'_J(h,k) = h'_J(h,kk') $$
$$( k ' \underset J \triangleleft (k \underset J \triangleright  h))(k \underset J \triangleleft h) = k'k \underset J \triangleleft h $$}
\begin{dm}
 The first identity is just lemma \ref{poubelles} where one flips $H$ and $K$. Also, for any $h \in H$ and $k,k' \in K$,  one has:

\begin{align*}
k'k \underset J \triangleleft  h
&= (k'k \underset J \triangleright  h)^{-1}k'kh = (k'k \underset J \triangleright  h)^{-1}k' (k \underset J \triangleright  h)(k \underset J \triangleleft  h)\\
&= (k'k \underset J \triangleright  h)^{-1}(k' \underset J \triangleright (k \underset J \triangleright  h))(k' \underset J \triangleleft (k \underset J \triangleright  h)(k \underset J \triangleleft  h)  \\
&= (k' \underset J \triangleleft (k \underset J \triangleright  h)(k \underset J \triangleleft  h)
\end{align*}
\end{dm}

\subsubsection{\bf{Theorem}}
\label{banlieuebis}
{\it   The family  $(W^J_k)_{k \in K}$ is a group of unitaries in $M'_1 \cap M_3$, it implements an action of $H$ on the sub $*$-algebra of $M'_1 \cap M_3$ generated by the family $(1 \otimes \lambda(s) \otimes \rho(s)\chi_h)_{s \in H \cap K, h\in K}$ which is equivalent to the action $\sigma^J$ defined in  \ref{domotermbis}. The $*$-algebra   $M'_1 \cap M_3$ is isomorphic to the crossed product $ ( C(H) \underset {\rho} \rtimes (H \cap K))\underset {\sigma^I}\rtimes K$ and to $\mathbb C  \mathcal T'$. } 
\vskip 0.2cm
\begin{dm}
For any $h_0 \in H$, and $k \in K$,  due to \ref{graphebis} ii), one has:
\begin{align*}
 W^J_k(1 \otimes 1 \otimes \chi_{h_0}) {W^J_k}^* 
 &=   \underset {h,h^1\in H } \sum v_{h,k,h'_J(h,k)}  (1 \otimes 1 \otimes \chi_{h_0})  v^*_{h^1,k,{h^1_J}'(h,k)}  
\\
 &=    v_{h_0,k,h'_J(h_0,k)}   v^*_{h_0,k,{h^1_J}'(h_0,k)}  \\
&=  1 \otimes 1 \otimes \chi_{ k \underset J \triangleright  h_0}
\end{align*}
Suming this equality for all $h_0 \in H$  gives that $W^J_k$ is a unitary.

\vskip 0.2cm

Now for any $s \in S$, $k \in K$ and $h \in H$, one easily sees that  $h'_J(hs^{-1},k)= sh'_J(h,k)s^{-1}$ and   $ k \underset J \triangleleft (hs^{-1})= s(k \underset J \triangleleft h)s^{-1}$, from this one  deduces that: 

\begin{align*}
(1 \otimes \lambda(s) &\otimes \rho(s))v_{h,k,h'_J(h,k)}
= (1 \otimes \lambda(s(k \underset J \triangleleft h)k^{-1}) \otimes \rho(sh'_J(h,k))\chi_{h}v_k\\
&= (1 \otimes \lambda(s(k \underset J \triangleleft h)s^{-1}k^{-1}ksk^{-1}) \otimes \rho(sh'_J(h,k)s^{-1}s)\chi_{h}v_k\\
&= (1 \otimes \lambda((k \underset J \triangleleft hs^{-1})k^{-1}ksk^{-1}) \otimes \rho(h'_J(hs^{-1},k)\rho(s))\chi_{h}v_k  \\
&=  (1 \otimes \lambda((k \underset J \triangleleft hs^{-1})k^{-1}) \otimes \rho(h'_J(hs^{-1},k)\chi_{hs^{-1}})( 1 \otimes \lambda(ksk^{-1}) \otimes \rho(s))v_k \\
&=  v_{hs^{-1},k,h'_J(hs^{-1},k)} \alpha_{-k}^2( 1 \otimes \lambda(ksk^{-1}) \otimes \rho(s))\\
&=  v_{hs^{-1},k,h'_J(hs^{-1},k)} ( 1 \otimes \lambda(s) \otimes \rho(s))
\end{align*}

From this one deduces that: 
\begin{align*}
W^J_k(1 \otimes \lambda(s) \otimes \rho(s)){W^J_k}^* = 1 \otimes \lambda(s) \otimes \rho(s)
 \end{align*}

\vskip 0.5cm

Now, for any $k,k'$ in $K$, due to \ref{graphebis} and \ref{bepc},  we have:

\begin{align*}
W^J_{k'}W^J_k
=     \underset{h,h_1 \in H } \sum     v_{h_1,k',h'_J(h_1,k')}v_{h,k,h'_J(h,k)}
=     \underset {h \in H} \sum  v_{ k \underset J \triangleright  h,k',h'_J(k \underset J \triangleright  h,k')}v_{h,k,h'_J(h,k)}
\end{align*}

But for any $h \in H$, one has:
\begin{align*}
&v_{ k \underset J \triangleright  h,k',h'_J(k \underset J \triangleright  h,k')}v_{h,k,h'_J(h,k)} = \\
&(1 \otimes \lambda((k' \underset J \triangleleft k \underset J \triangleright  h)k'^{-1}) \otimes \rho(h'_J(k \underset J \triangleright  h,k'))\chi_{k \underset J \triangleright  h})v_{k'}(1 \otimes \lambda((k \underset J \triangleleft h)k^{-1}) \otimes \rho(h'_J(h,k))\chi_{h})v_k  \\
&= (1 \otimes \lambda((k' \underset J \triangleleft k \underset J \triangleright  h)k'^{-1}) \otimes \rho(h'_J(k \underset J \triangleright  h,k')))v_{k'}(1 \otimes \lambda((k \underset J \triangleleft h)k^{-1}) \otimes \rho(h'_J(h,k))\chi_{h})v_k  \\
&= (1 \otimes \lambda((k' \underset J \triangleleft k \underset J \triangleright  h)k'^{-1}) \otimes \rho(h'_J(k \underset J \triangleright  h,k')))(1 \otimes \lambda(k'(k \underset J \triangleleft h)({k'k})^{-1}) \otimes \rho(h'_J(h,k))\chi_{h})v_{kk'}  \\
&= (1 \otimes \lambda((k' \underset J \triangleleft k \underset J \triangleright  h)(k \underset J \triangleleft h)({k'k})^{-1} )\otimes \rho(h'_J(k \underset J \triangleright  h,k')h'_J(h,k))\chi_{h})v_{kk'} \\
&= (1 \otimes \lambda((k' k \underset J \triangleleft  h)({k'k})^{-1} )\otimes \rho(h'_J(h,k'k))\chi_{h})v_{kk'} 
\end{align*}

One deduces that:

\begin{align*}
W^J_{k'}W^J_k
=  \underset{h \in H } \sum (1 \otimes \lambda((k' k \underset J \triangleleft  h)({k'k})^{-1} )\otimes \rho(h'_J(h,k'k))\chi_{h})v_{kk'}  = W^J_{k'k}
\end{align*}

So  the family $(W^J_k)_{k \in K}$ is a group of unitaries, which implements an action of $H$ on the sub $*$-algebra of $M'_1 \cap M_3$ generated by the family $(1 \otimes \lambda(s) \otimes \rho(s)\chi_h)_{s \in S, h\in K}$,  as for any $s  \in S$, $h \in H$ and $kK$, one has: $(1 \otimes \lambda (s)\otimes  \rho(s) )W^J_k( 1 \otimes 1 \otimes  \chi_{h}) = v_{h ,k, sh'_J(h,k)}$, this proves that the crossed product  is exactly $M'_1 \cap M_3$ and so, by theorem \ref{roy} where $H$ and $K$ are flipped, $M'_1 \cap M_3$ is isomorphic to the crossed product $ ( C(K) \underset {\rho} \rtimes S)\underset {\sigma^I_h}\rtimes H $. \end{dm}

\vskip 1cm
  \subsubsection{\bf {Corollary and notations}}
\label{jesouffrebis}
{\it For any $t \in \mathcal T'$, let $(k,h,s)$ be the unique element 
\vskip 0,3cm
in $K \times H\times S$ such that $t = \  \    k \underset J\triangleright h   \eenmatrix{k}{hs}{ }{(k \underset J \triangleleft  h)s  } \   \ $, let $\theta^J_t$ be the element in $M'_1 \cap M_3$
\vskip 0,2cm
 equal to $W^J_k(1 \otimes \lambda(s) \otimes \chi_h\rho(s) )$, then the family $(\theta^J_t)_{t \in \mathcal T}$ is a basis of $M'_1 \cap M_3$ such that for any $t' \in \mathcal T'$, one has:

$$\theta^J_t\theta^J_{t'} =  \   \ \left \{ \begin{array}{rl}
& \theta^J_{ t \underset H{\overset h \star}t'} \  \mathrm{if \ \  t \ \ and\ \ t' \ \ are \ \ composable\ \  for} \ \ \underset K{\overset h \star}  \hskip 0.5cm   \\
&  0 \  \ \mathrm{otherwise} \hskip 0.5cm  
\end{array}
\right. $$
$$ ( \theta_t)^* = \theta_{t^{-h}} $$
}
\newline
\begin{dm}
Using theorem \ref{banlieuebis}, this is just a reformulation of proposition \ref{calmos}
\end{dm}

\vskip 2cm
\subsection{The  co-algebras structures  of $M'_0 \cap M_2$ and $M'_1 \cap M_3$}

Let's apply the results given in paragraph \ref{ses} to find  co-algebras structures  on $M'_0 \cap M_2$ and $M'_1 \cap M_3$. With notations of  paragraph \ref{temoin}, one can use  the  pairing with $M'_1 \cap M_3$ defined for any $a \in M'_0 \cap M_2$ and $b \in M'_1 \cap M_3$  by: 
$$<a,b> = |H \cap K||H|^2|K|^2\tau(ae_2e_1b).$$
For $M'_0 \cap M_2$, the coproduct $\Gamma$, the antipod $\kappa$ and the counit $\epsilon$ are given by the following formulas:
\begin{align*}
 & \hskip 0.5cm \epsilon(a) = <a,1>   \\  
 &<\Gamma(a), b \otimes b'> = <a,bb'> \\
& <\kappa(a), b> = \overline{<a^*,b^*>}\end{align*}

One has equivalent formulas for $(M'_1 \cap M_3, \hat \Gamma,\hat \kappa,\hat \epsilon)$.

\subsubsection{\bf{Remark}}
The general bracket given in \cite{NV1} or \cite{D}, uses a slightly more complicated formula, as we shall see later we are here in a simplier situation for which $\Gamma$ is multiplicative and    one can apply theorem 4.17 of \cite{NV1}.
\vskip 1cm

\subsubsection{{\bf Lemma}}
\label{espagne}
{\it  For any $h,h'\in H$, $k,\in K$, $s,s'\in S$, one has:
\begin{align*}
(1 \otimes \rho(s) &\chi_k \otimes 1)W^I_he_2e_1(1 \otimes \lambda(s') \otimes \rho(s')\chi_{h'}) = \\
&= |H|^{-1}|K|^{-1}v_{h^{-1}  \underset I \triangleright'k}(1 \otimes \rho(sk^{-1}(h^{-1}  \underset I \triangleright' k ))\chi_e \lambda(s') \otimes \rho(hh')\chi_{h'})
\end{align*} }
\vskip 0.7cm
\begin{dm}
 For any $h\in H$, $k,\in K$, $s\in S$, one has:
 \begin{align*}
 &(1 \otimes \rho(s)   \chi_k \otimes 1)W^I_he_2e_1=\\
  &= |K|^{-1}(1 \otimes \rho(s) \otimes 1)(1 \otimes \rho(k^{-1}(h^{-1}  \underset I \triangleright' k ))\chi_{h^{-1}  \underset I \triangleright' k }\otimes \rho(h))(1 \otimes 1 \otimes \chi_e)( \underset {k_1 \in K} \sum v_{k_1})e_1\\
  &= |H|^{-1}|K|^{-1} \underset {k_1 \in K} \sum v_{k_1} \alpha^2_{k_1^{-1}}(1 \otimes \rho(sk^{-1}(h^{-1}  \underset I \triangleright' k ))\chi_{h^{-1}  \underset I \triangleright' k }\otimes \rho(h) \chi_e)(1 \otimes \chi_e \otimes \underset {h_1 \in H} \sum \rho(h_1))\\
  &= |H|^{-1}|K|^{-1} \underset {k_1 \in K,h_1\in H} \sum v_{k_1} (1 \otimes \rho(sk^{-1}(h^{-1}  \underset I \triangleright' k ))\chi_{k_1^{-1}(h^{-1}  \underset I \triangleright' k) }\chi_e \otimes \rho(h) \chi_e   \rho(h_1))\\  &= |H|^{-1}|K|^{-1}v_{h^{-1}  \underset I \triangleright' k}(1 \otimes \rho(sk^{-1}(h^{-1}  \underset I \triangleright'k ))\chi_e \otimes \rho(h)\chi_e(\underset {h_1\in H} \sum \rho(h_1)))
 \end{align*}
 
 Hence, for any $h' \in H$ and $s' \in S$, one has:
 
 \begin{align*}
 &(1 \otimes \rho(s)   \chi_k \otimes 1)W^I_he_2e_1(1 \otimes \lambda(s') \otimes \rho(s') \chi_{h'})= \\
   &= |H|^{-1}|K|^{-1}v_{h^{-1}  \underset I \triangleright' k}(1 \otimes \rho(sk^{-1}(h^{-1}  \underset I \triangleright' k ))\chi_e\lambda(s')  \otimes \rho(h)\chi_e(\underset {h_1\in H} \sum \rho(h_1)) \rho(s')\chi_{h'})\\
   &= |H|^{-1}|K|^{-1}v_{h^{-1}  \underset I \triangleleft k}(1 \otimes \rho(sk^{-1}(h^{-1}  \underset I \triangleright'k ))\chi_e\lambda(s')  \otimes \rho(h)\chi_e(\underset {h_1\in H} \sum \rho(h_1s'))\chi_{h'})\\
   &= |H|^{-1}|K|^{-1}v_{h^{-1}  \underset I \triangleright' k}(1 \otimes \rho(sk^{-1}(h^{-1}  \underset I \triangleright' k ))\chi_e\lambda(s')  \otimes \rho(h)\chi_e(\underset {h_2\in H} \sum \rho(h_2))\chi_{h'})\\
   &= |H|^{-1}|K|^{-1}v_{h^{-1}  \underset I \triangleright' k}(1 \otimes \rho(sk^{-1}(h^{-1}  \underset I \triangleright' k ))\chi_e\lambda(s')  \otimes \rho(h)(\underset {h_2\in H} \sum \rho(h_2)\chi_{h_2}\chi_{h'})\\
   & = |H|^{-1}|K|^{-1}v_{h^{-1}  \underset I \triangleright' k}(1 \otimes \rho(sk^{-1}(h^{-1}  \underset I \triangleright' k ))\chi_e \lambda(s') \otimes \rho(hh')\chi_{h'})
    \end{align*}
 
\end{dm}

\vskip 1cm

\subsubsection{{\bf Lemma}}
\label{portugal}
{\it  For any $h,h' \in H$, $k,k' \in K$, $s,s' \in S$, one has:
\begin{align*}
  <(1 \otimes \rho(s)\chi_k \otimes 1&)W^I_h,(1 \otimes \lambda(s') \otimes \rho(s')  \chi_{h'}) W^J_{k'} >= \\
 & = |S| \delta_{k'h^{-1}  \underset I \triangleright' k,e}\delta_{h(k'^{-1} \underset J \triangleright h'),e}\delta_{(k' \underset J \triangleleft (k'^{-1} \underset J \triangleright  h'))k,{s'}^{-1}s}   
  \end{align*} }
\vskip 0.7cm
\begin{dm}
 For any $h' \in H$, $k' \in K$, due to \ref{domotermbis} and  \ref{banlieue}, one has:
\begin{align*}
(1\otimes 1 \otimes \delta_{h'} )  W^J_{k'} &=  (1\otimes 1 \otimes \delta_{h'} ) W^J_{k'} ( 1 \otimes 1 \otimes  \delta_{k'^{-1} \underset J \triangleright h'} )\\
&= \ (1\otimes 1 \otimes \delta_{h'} )(1 \otimes \Theta_1 \otimes \Theta_2 \chi_{k'^{-1} \underset J \triangleright h'})v_{k'} 
\end{align*}

where $\Theta_1 = \lambda((k' \underset J \triangleleft (k'^{-1} \underset J \triangleright  h'))k'^{-1})$ and $\Theta_2 = \rho(h'^{-1}(k'^{-1} \underset J \triangleright h'))$.

Hence, using lemma \ref{espagne}, one has:

\begin{align*}
&(1 \otimes \rho(s)   \chi_k \otimes 1)W^I_he_2e_1(1 \otimes \lambda(s') \otimes \rho(s') \chi_{h'})   W^J_{k'} =\\
&= |H|^{-1}|K|^{-1}v_{h^{-1}  \underset I \triangleright' k}(1 \otimes \rho(sk^{-1}(h^{-1}  \underset I \triangleright' k ))\chi_e \lambda(s')\Theta_1 \otimes \rho(hh')\chi_{h'}\Theta_2\chi_{k'^{-1} \underset J \triangleright h'})v_{k'}  \\
&= |H|^{-1}|K|^{-1}v_{h^{-1}  \underset I \triangleright' k}(1 \otimes \rho(sk^{-1}(h^{-1}  \underset I \triangleright' k ))\chi_e \lambda(s')\Theta_1 \otimes \rho(h(k'^{-1} \underset J \triangleright h'))\chi_{k'^{-1} \underset J \triangleright h'})v_{k'} 
\end{align*}

As $\tau$ is a trace, one has:
\begin{align*}
&(|H||K||H\cap K|)^{-1}<(1 \otimes \rho(s)\chi_k \otimes 1)W^I_h,(1 \otimes \lambda(s') \otimes \rho(s')  \chi_{h'})  W^J_{k'} >=\\
&= \tau((1 \otimes \rho(sk^{-1}(h^{-1}  \underset I \triangleright' k ))\chi_e \lambda(s')\Theta_1 \otimes \rho(h(k'^{-1} \underset J \triangleright h'))\chi_{k'^{-1} \underset J \triangleright h'})v_{k'}v_{h^{-1}  \underset I \triangleright'k})\\
&=\tau(\beta^1(1 \otimes \rho(sk^{-1}(h^{-1}  \underset I \triangleright' k ))\chi_e \lambda(s')\Theta_1 )(1 \otimes 1 \otimes \rho(h(k'^{-1} \underset J \triangleright h'))\chi_{k'^{-1} \underset J \triangleright h'})v_{k'h^{-1}  \underset I \triangleleft k})\\
&= |K|^{-1}\delta_{k'h^{-1}  \underset I \triangleright' k,e}\delta_{h(k'^{-1} \underset J \triangleright h'),e}\tau(\beta^1(1 \otimes \rho(sk^{-1}(h^{-1}  \underset I \triangleright' k ))\chi_e \lambda(s')\Theta_1 ))
\end{align*}
\vskip 0.4cm
Using the fact that for any $k_1 \in K$ $\chi_{e} \lambda(k_1) =  \chi_{e}\rho(k_1^{-1}) = \rho(k_1^{-1})\chi_{k_1^{-1}}$, one also has:

\begin{align*}
&(|H||K||H\cap K|)^{-1}<(1 \otimes \rho(s)\chi_k \otimes 1)W^I_h,(1 \otimes \lambda(s') \otimes \rho(s')  \chi_{h'}) W^J_{k'} >=\\
&= |K|^{-1}\delta_{k'h^{-1}  \underset I \triangleright' k,e}\delta_{h(k'^{-1} \underset J \triangleright h'),e}\tau(\beta^1(1 \otimes \rho(sk^{-1}(h^{-1}  \underset I \triangleright' k )k'(k' \underset J \triangleleft (k'^{-1} \underset J \triangleright  h'))^{-1}{s'}^{-1}) \chi_{z}))
\end{align*}
for a certain $z$, so we have:
\begin{align*}
&(|H||K||H\cap K|)^{-1}<(1 \otimes \rho(s)\chi_k \otimes 1)W^I_h,(1 \otimes \lambda(s') \otimes \rho(s')  \chi_{h'})  W^J_{k'} >=\\
&= |K|^{-1}\delta_{k'h^{-1}  \underset I \triangleright' k,e}\delta_{h(k'^{-1} \underset J \triangleright h'),e}\tau(\beta^1(1 \otimes \rho(sk^{-1}(k' \underset J \triangleleft (k'^{-1} \underset J \triangleright  h'))^{-1}{s'}^{-1}) \chi_{z}))\\
&= |K|^{-1}|H|^{-1}\delta_{k'h^{-1}  \underset I \triangleright' k,e}\delta_{h(k'^{-1} \underset J \triangleright h'),e}\delta_{sk^{-1}(k' \underset J \triangleleft (k'^{-1} \underset J \triangleright  h'))^{-1}{s'}^{-1},e}\\
&= |K|^{-1}|H|^{-1}\delta_{k'h^{-1}  \underset I \triangleright' k,e}\delta_{h(k'^{-1} \underset J \triangleright h'),e}\delta_{(k' \underset J \triangleleft (k'^{-1} \underset J \triangleright  h'))k,{s'}^{-1}s}
\end{align*} 
The lemma follows
\end{dm}
\vskip 1.1cm

\subsubsection{{\bf Lemma}}
\label{deprime} 
{\it    For all $h,h' \in H$ and $k,k' \in K$ such that  $k'= h   \underset I \triangleright' k$ and $h= k' \underset J \triangleright h'$, then there exists a unique $\sigma \in S$ such that: $ \sigma  = (k'h')^{-1}hk = ({k'}^{-1} \underset J \triangleleft (k'  \underset J \triangleright  h'))k$ }
\vskip 0,6cm
\begin{dm}
For all $h,h' \in H$ and $k,k' \in K$ such that  $k'= h   \underset I \triangleright' k$ and $h= k' \underset J \triangleright h'$, then there exist $h_1 \in H$ and $k_1\in K$ such that $hk_1= k'h'$ and $hk = k'h_1$. One deduces that $h'^{-1}k'^{-1}hk = k_1^{-1}k = h'^{-1}h_1$ and the existence (and uniqueness)  of $\sigma$ such that $ \sigma  = (k'h')^{-1}hk$. One also has:
\begin{align*}
({k'}^{-1} \underset J \triangleleft (k'  \underset J \triangleright  h'))k
&= ({k'}^{-1} \underset J \triangleleft h)k = ({k'}^{-1} \underset J \triangleright h)^{-1}({k'}^{-1} h)k =  h'^{-1}({k'}^{-1} h)k = (k'h')^{-1}hk
\end{align*}
The lemma follows
\end{dm}
 \vskip 1.5cm

\subsubsection{{\bf Proposition}}
\label{artifice}
 { \it Using notations \ref{inverse}, \ref{jesouffre} and \ref{jesouffrebis}, for any $x \in \mathcal T$ and  $x' \in \mathcal T'$, one has:
$$
  <  \theta^I_x ,  \theta^J_{x'}> = |H \cap K|\delta_{x^t,x'} = <x,x'>
$$}

\begin{dm}
For any $x \in \mathcal T$ and  $x' \in \mathcal T'$ let $ (h,k,s)$ in $H\times   K \times S$ and $ (k', h', s')$ 
\vskip 0.3cm 
in $K \times H \times S$ such that $x = \  \    h \underset I\triangleright' k   \eenmatrix{h}{ks}{ }{(h \underset I \triangleleft'  k)s  } \   \ $ and $x' = k' \underset J\triangleright h'   \eenmatrix{k'}{h's'}{ }{(k' \underset J \triangleleft'  h')s'  }$, \  \  , then due to lemmas \ref{portugal}, \ref{deprime} and \ref{lafinarrivebis} one has:
\vskip 0.3cm 
\begin{align*}
    < \theta^I_{x},   \theta^J_{x'} >
  &=     |S| \delta_{k',h   \underset I \triangleright' k}\delta_{h, k'  \underset J \triangleright h'}\delta_{({k'}^{-1} \underset J \triangleleft (k'  \underset J \triangleright  h'))k,s'{s}^{-1}} \\
  &= |S| \delta_{k',h   \underset I \triangleright' k}\delta_{h, k'  \underset J \triangleright h'}\delta_{(k'h')^{-1}hk,s'{s}^{-1}}\\
  &= |S| \delta_{k',h   \underset I \triangleright' k}\delta_{h, k'  \underset J \triangleright h'}\delta_{hks,k'h's'}\\
   &= \left \{ \begin{array}{rl}
& |S | \hskip 0.5cm  \mathrm{if} \   \left \{ \begin{array}{rl}
&hks = k'h's' \hskip 0.5cm   \\
&k' = h \underset I \triangleright' k  \hskip 0.5cm \\
&h = k' \underset J \triangleright h' \hskip 0.5cm 
\end{array}
\right. \\
&0 \hskip 0.5cm  \ \mathrm{otherwise}
\end{array}
\right.\\
&= |S|\delta_{x^t,x'} \\
&= <x,x'>
\end{align*}
\end{dm}

\subsubsection{{\bf Theorem}}
\label{train} 
{\it The pair of C*-quantum groupoids in duality  $(M'_0 \cap M_2,\Gamma,\kappa,\epsilon)$ and $(M'_1 \cap M_3,\hat \Gamma,\hat \kappa,\hat \epsilon)$ is isomorphic to the pair of C*-quantum groupoids in duality associated with $(\mathbb C \mathcal T, \mathbb C \mathcal T')$. For any $ x \in \mathcal T$, and $x' \in \mathcal T'$, one has:
\begin{align*}
\Gamma(\theta^I_{x})  =  \frac{1}{| H \cap K |} \underset{ x_2\underset H{\overset {h} \star} x_1 = x} \sum \theta^I_{x_1} \otimes \theta^I_{x_2} \ \ ; \ \ \hat \Gamma(\theta^J_{x})  =  \frac{1}{| H \cap K |} \underset{ x_2\underset K{\overset {h} \star} x_1 = x'} \sum \theta^J_{x_1} \otimes \theta^J_{x_2}
\end{align*}
\begin{align*}
\kappa(\theta^I_{x}) = \theta^I_{ x^{-hv}} \ \ ; \ \ \hat\kappa(\theta^J_{x'}) = \theta^J_{ x'^{-hv}} \end{align*}  \begin{align*}
&\epsilon(\theta^I_{x}) = \left \{ \begin{array}{rl}
& |H \cap K | \hskip 0.5cm  \mathrm{if \  \ x  \ \ of\ \ the\ \ form } \ \  \eenmatrixbase{h}{e}{ e}{h   } \\
&0 \hskip 0.5cm  \ \mathrm{otherwise}
\end{array}
\right.
\end{align*}
\begin{align*}
&\hat \epsilon(\theta^J_{x}) = \left \{ \begin{array}{rl}
& |H \cap K | \hskip 0.5cm  \mathrm{if \  \ x'  \ \ of\ \ the\ \ form } \ \  \eenmatrixbase{k}{e}{ e}{k   } \\
&0 \hskip 0.5cm  \ \mathrm{otherwise}
\end{array}
\right.
\end{align*}
\begin{align*}
 <\theta^I_x,\theta^J_{x'}> = |H \cap K|\delta_{x^t,x'} = <x,x'>
\end{align*}
}
\begin{dm}  
Obvious by Proposition \ref{artifice}
 \end{dm}

\subsubsection{{\bf Corollary}}
\label{fin} 
{\it The von Neumann algebra crossed product $\mathbb C \mathcal T \rtimes \mathbb C \mathcal T'$ is isomorphic to $\mathbb C[H \cap K] \otimes \mathcal L(\mathbb C^{|H||K|})$.}
\vskip 0.5cm
\begin{dm}
 This is a direct consequence of corollary \ref{jncg} and corollary \ref{indiana}.
\end{dm}


\clearpage

\addcontentsline{toc}{chapter}{Bibliographie}

\begin{thebibliography}{30}
\bibitem[AN1]{AN1} N.Andruskiewitsch \& S.Natale : Double categories and  quantum  groupoids {\it Publ. Mat. Urug. } {\bf 10}  (2005) 11-51;

\bibitem[AN2]{AN2} N.Andruskiewitsch \& S.Natale : Tensor categories attached to double groupoids  (math Q.A. 0408045) to appear in Adv. Math  ;

\bibitem[B] {B} E.J.Beggs : Making non-trivially associated tensor categories from left coset representatives { \it Journal of pure and Applied Algebra} {\bf 177} (2003), 5-41;

\bibitem[B1]{B1}  D.Bisch :  Higher Relative Commutants and the Fusion Algebra Associated to a Subfactor { \it Fields Institute Communications} {\bf 13} (1997), 13-63 ;

\bibitem[BS]{BS} S. Baaj \& G. Skandalis :Unitaires multiplicatifs et dualit\'e pour 
les produits crois\'es de C*-alg\`ebres. {\it Ann. Sci. ENS} {\bf 26} 
(1993), 425-488;

\bibitem[BBS]{BBS} S. Baaj  \& E. Blanchard \& G. Skandalis :Unitaires multiplicatifs 
en dimension finie et
leurs sous-objets. {\it Ann.Inst. Fourier} {\bf 49} (1999), 1305-1344 ;

\bibitem[BH]{BH} D.Bisch \& U.Haagerup :  Composition of subfactors: new examples of infinite depth subfactors, {\it Ann. Sci. ENS 4\`eme  s\' erie } {\bf 29} $n^o 3$ (1996), 329-383 ;

\bibitem[BoSz]{BoSz} G. B\"ohm \& K.Szlach\'anyi :  Weak C*-Hopf algebras: the  
coassociative symmetry of non integral dimensions, {\it  Quantum groups 
and quantum spaces. Banach Center Publications} {\bf 40} (1997), 9-19  ;

\bibitem[BoSzNi] {BoSzNi} G. B\"ohm \& K.Szlach\'anyi \& F.Nill : Weak Hopf Algebras I. 
Integral Theory and $C^*$-structure. {\it Journal of Algebra} {\bf 221} 
(1999),  385-438  ;

\bibitem[D]{D} M.C. David : $C^*$-groupoides quantiques et inclusions de facteurs: structure symŽtrique et autodualitŽ actions sur le facteur hyperfini de type $II_1$. {\it Journal of operator theory} {\bf 54:1}(2005), 27-68;

\bibitem[D2]{D2} M.C. David : Private communication

\bibitem[E1]{E1} M. Enock :  Inclusions of von Neumann algebras and quantum groupoids III .
{\it Journal of Functional Analysis} {\bf 223} (2005), 311-364  ;

\bibitem[E2]{E2} M. Enock :  Measured quatum groupoids in action ;

\bibitem[EH]{EH} C.Ehresmann:  Cat\'egories structur\'ees. {\it  Ann.Sci.Ecole Norm.Sup.} (3), {\bf 80} (1963) 349-426 ;

\bibitem[EV]{EV} M. Enock  \& J.M. Vallin : Inclusions of von Neumann
algebras and quantum groupoids,
{\it Journal of Functional Analysis {\bf 172} (2000), 249-300.} ;


\bibitem[G]{G} C.F.Gardiner, Algebraic structures : {\it Ellis Horwood limited} ({\bf1986}) John Wiley \& Sons ;

\bibitem[GHJ]{GHJ} F.M.Goodman, P. de la Harpe, V.F.R. Jones :  Coxeter graphs and towers of algebras {\it M.S.R.I. publications} 14;
   
\bibitem[HS] {HS}J.H.Hong , W.Szymanski :  Composition of subfactors and twisted crossed products.{\it Journal of operator theory} {\bf 37}(1997), 281-302. ;


\bibitem[L]{L} F.Lesieur:  thesis, http://tel.ccsd.cnrs.fr/documents/archives0/00/00/55/05  ;


\bibitem[NV1]{NV1} D. Nikshysh \& L. Vainerman :  A characterization of depth 2 
subfactors of $II_1$ factors. {\it JFA} {\bf 171} (2000), 278-307  ;

\bibitem[NV2]{NV2} D. Nikshysh \& L. Vainerman :  A Galois correspondence for 
$II_1$-factors and quantum groupoids. {\it JFA} {\bf 178} (2000),
113-142 ;

\bibitem[PP]{PP} M.Pimsner  \& S.Popa :  Iterating the basic construction, Trans.Amer.Math.Soc. {\bf 310}(1988), 127-133 ;


\bibitem[R]{R} J. Renault :  {\it A groupoid approach to $C^*$-algebras}   Lect.Notes
in Mah. {\bf 793} Springer-Verlag 1980 ;

\bibitem[Val1]{Val1} J.M. Vallin :  Groupo\" \i des quantiques finis. {\it Journal 
of Algebra} {\bf 26} (2001), 425-488  ;

\bibitem[Val2]{Val2}  J.M. Vallin :   Actions and coactions of finite quantum groupoids on von Neumann algebras, extensions of the  matched pair procedure ( submitted to Journal of Algebra)  ;

\bibitem[Val3]{Val3} J.M. Vallin :  $C^*$-alg\`ebres de Hopf et $C^*$-alg\`ebres de Kac. {\it Pro.London.Math.Soc} {\bf 50}, No.3 (1985), 131-174 ;

\bibitem[VV] {VV}S.Vaes \& L.Vainerman :  Extensions of locally compact quantum groups and the bicrossed product construction { \it Advances in Mathematics} {\bf 175} (2003), 1-101  ; \end{thebibliography}
\end{document}